\documentclass{amsart}
\usepackage{amsmath,amssymb}

\title[Geometrically infinite groups]{Constructing geometrically infinite groups on boundaries of 
deformation spaces}
\newtheorem{lem}{Lemma}[section] 
\newtheorem{thm}[lem]{Theorem}
\newtheorem{cor}[lem]{Corollary} 
\newtheorem{prop}[lem]{Proposition}
\newtheorem{ththm}{Theorem}

\theoremstyle{definition}
\newtheorem{dfn}[lem]{Definition}
\newtheorem{rem}[lem]{Remark}

\author{Ken'ichi Ohshika}
\address{Department of Mathematics, Graduate School of Science, 
Osaka University, Toyonaka, Osaka 560-0043, Japan}
\email{ohshika@math.sci.osaka-u.ac.jp}
\def\C{\mathbf{C}}
\def\H{\mathbf{H}}
\def\R{\mathbf{R}}

\def\ML{\mathcal{ML}}

\def\l{\mathrm{length}}
\def\T{\mathcal{T}}
\def\I{\mathrm{Int}}
\def\G{\Gamma}
\def\M{\mathcal{M}}
\def\PL{\mathcal{PML}}
\def\s{\mathcal{S}}

\def\Z{\mathbf{Z}}

\keywords{Kleinian group, deformation space, geometrically finite group}
\begin{document}
\maketitle

\begin{abstract}
Consider a geometrically finite Kleinian group $G$ without  parabolic or elliptic elements, with its Kleinian manifold
$M=(\H^3\cup \Omega_G)/G$.
Suppose that for each boundary component of $M$, either a maximal and connected
measured
lamination in the Masur domain or a marked conformal structure is given.
In this setting, we shall prove that there is an algebraic limit $\Gamma$ of
quasi-conformal deformations of $G$ such that there is a homeomorphism $h$ from $\mathrm{Int} M$ to $\H^3/\Gamma$ compatible with the natural isomorphism from $G$ to $\Gamma$,  the given laminations are unrealisable in $\H^3/\Gamma$, and the given conformal structures are pushed forward by $h$ to those  of $\H^3/\Gamma$.
Based on this theorem and its proof,  in the subsequent paper, the Bers-Thurston
conjecture, saying that every finitely generated
Kleinian group is an algebraic limit of quasi-conformal deformations of 
 minimally parabolic geometrically finite group, is proved using recent
solutions of Marden's conjecture by Agol, Calegari-Gabai, and the ending
lamination conjecture by Minsky collaborating with Brock, Canary and Masur.
\end{abstract}
\section{Introduction}
The deformation space $AH(G)$ of a Kleinian group $G$ is defined to be 
the space of faithful discrete representations of $G$ into $PSL_2\C$ 
taking parabolic elements to parabolic elements, modulo conjugacy.
One of the most important topics in the theory of Kleinian groups is 
to study topological structure of deformation spaces.
When a Kleinian group $G$ is geometrically finite, it is known that 
its quasi-conformal deformation space, which we denote by $QH(G)$, 
forms an open neighbourhood of the identity representation of $G$ in its entire deformation
space 
$AH(G)$.
The topological types of quasi-conformal deformation spaces are 
well understood, in terms of Teichm\"{u}ller space of the region of 
discontinuity,  by work of Ahlfors, Bers, Kra, Maskit and Sullivan 
among others.
It is conjectured by Bers and Thurston that the entire deformation 
space $AH(G)$ is the closure of the union of quasi-conformal 
deformation spaces of (minimally parabolic) geometrically finite Kleinian groups
isomorphic 
to $G$.
Taking this into account, we can see that a crucial step for 
understanding deformation spaces is to study the boundaries of 
quasi-conformal deformation spaces.

On the other hand, the ending lamination conjecture, which was recently proved by Minsky  partially in collaboration  with Brock, Canary and Masur (\cite{Mi}, \cite{BCM}), asserts that topologically tame Kleinian groups can be classified, up to conjugacy, completely by  the following four invariants: the homeomorphism types of the  quotient hyperbolic 3-manifolds, the parabolic locus, the conformal structures at infinity  associated to geometrically finite ends, and the ending laminations associated to geometrically infinite ends.
Also by recent resolution of Marden's conjecture by Agol and Calegari-Gabai, it is now known that every finitely generated Kleinian group is topologically tame.
Consider a geometrically finite Kleinian group $G$ and its quasi-conformal deformation space $QH(G)$.
Then, to obtain a complete list of topologically tame Kleinian groups in $AH(G)$, we have only to show that for a given quadruple of  invariants, with obvious necessary conditions for the quadruple to be realised, such as the condition that no two laminations are homotopic in the $3$-manifold, we  can construct a (topologically tame) Kleinian group on the boundary of $QH(G)$ realising the triple.

In the main theorem of this paper, we shall prove that if we assume that the parabolic locus is empty, any triple of the remaining invariants, such that the third one consists of maximal and connected laminations represented by a measured lamination in the Masur domain no two of which are homotopic and none of which cover laminations on embedded non-orientable surfaces, can be realised by a topologically tame Kleinian group without
 parabolic elements on the boundary of the quasi-conformal deformation space as its homeomorphism type, its conformal structures at infinity, and its unrealisable laminations. 
 
We should note the following facts to understand the meaning of this statement.
As was shown in \cite{OhM}, the first invariant always coincides with 
the homeomorphism type of $\H^3/G$ for a topologically 
tame Kleinian 
group without parabolic elements on the boundary of the 
quasi-conformal deformation space.
(This is not the case if we allow new parabolic elements to appear in the  limit.
See Anderson-Canary \cite{AC}.)
Therefore, we need to fix the first invariant to be the homeomorphism type of
$\H^3/G$.
By a result of Canary \cite{CaJ}, for purely loxodromic topologically tame
groups, every
measured lamination representing an
ending lamination must be a maximal and connected lamination contained in the 
Masur domain, which is the set of measured laminations intersecting every limit of boundaries of
compressing discs. 

Two ending laminations for distinct ends cannot be 
homotopic each other and no ending lamination can be homotopic to a double cover
of a
lamination on a closed non-orientable surface,
as can be seen by the argument of
\S 6.4 in 
Bonahon \cite{BoA}.
In our main theorem, it will turn out that  these conditions for measured laminations on the boundary, that every one is contained in the Masur domain, that no two are homotopic, and that no one doubly covers a lamination on a non-orientable surface, are 
also sufficient for them to represent unrealisable laminations in a topologically tame 
Kleinian group on the boundary of the quasi-conformal deformation space.

By Bonahon's work,  ending laminations must be unrealisable laminations.
Conversely, it can be proved that  any unrealisable lamination contained in the Masur domain of a boundary component of a core  is the image of an ending lamination by an auto-homeomorphism acting on the fundamental group by an inner-automorphism.
This fact, which should be regarded as the uniqueness of unrealisable laminations, was proved by Bonahon in the case of freely indecomposable groups.
For non-free but freely decomposable groups the argument in Ohshika \cite{OhM} shows this.
It is most difficult to prove this fact in the case when the group is free.
In the setting of the limit group as in our main theorem, this has been proved in Ohshika \cite{Ohpr}, which is the sequel to the present paper.
For free Kleinian groups in general, Namazi and Souto \cite{NS} have proved this.

Based on the argument in this paper,  in \cite{Ohpr}, we have succeeded in
generalising the
convergence theorem which is the key step of our argument in this paper by
removing
the restrictions that $G$ is purely loxodromic and that given laminations are
maximal and connected.
Combining this result with the resolutions of Marden's conjecture and the ending
lamination conjecture, and showing that the unrealised laminations in the
limit groups are actually ending laminations as was explained above, we have proved that the
Bers-Thurston density conjecture is true.
Namely, {\sl every finitely generated Kleinian group is an algebraic limit of
quasi-conformal deformations of a minimally parabolic geometrically finite
group.}
See Ohshika \cite{Ohpr} for more details.

For the Bers-Thurston density conjecture in the case of freely indecomposable groups, there is another approach which is quite different from ours.
This is due to Bromberg and Bromberg-Brock (\cite{Brom}, \cite{BB}) and used the technique of the deformation of cone manifolds.
Their method uses the ending lamination conjecture only for the case when the manifolds have injectivity radii bounded away from $0$.

Our main theorem can also be regarded as a generalisation of the main theorem in 
\cite{OhL}, where we dealt with Kleinian groups whose corresponding 
hyperbolic 3-manifolds have boundary-irreducible compact cores.
To make the result of \cite{OhL} adaptable to general (finitely generated and 
torsion-free) Kleinian groups, it was 
necessary to prove a convergence theorem for freely decomposable 
Kleinian groups, generalising the main theorem of \cite{OhI}.
One of the most difficult steps for that had been the case of functions groups, 
which was solved by Kleineidam-Souto \cite{KS} recently.
See Theorem \ref{KS}.
 In \cite{KS2}, they also gave a proof  of a special case of our main theorem in
the case when $G$ is a geometrically
finite function group independently of our work here.

Now we state the main theorem.

\begin{thm}
\label{main}
Let $G$ be a geometrically finite Kleinian group without torsions 
or parabolic elements.
Let $M$ be a compact 3-manifold whose interior is homeomorphic to 
$\mathbf{H}^3/G$, and $S_1, \ldots , S_m$ its boundary components.
(We can identify $M$ with the Kleinian manifold $(\H^3 \cup 
\Omega_G)/G$.)
Suppose that at least one maximal and connected measured laminations $\lambda_{j_1}, \ldots , 
\lambda_{j_p}$ contained in the Masur domains are given on  boundary components
$S_{j_1} ,
\ldots , 
S_{j_p}$ among $S_1,\ldots S_m$, and marked conformal structures $m_1, 
\ldots , m_q$ on the 
remaining boundary components $S_{i_1} , \ldots , S_{i_q}$. 
When $M$ is homeomorphic to $S_{j_1} \times I$ and $p=2$, we further 
assume that the supports of  $\lambda_{j_1}$ and $\lambda_{j_2}$ are 
not homotopic in $M$.
When  $M$ is  homeomorphic to a twisted $I$-bundle over a non-orientable
surface $S'$,
we further assume that $\lambda_{j_1}$ is not a lift of a measured lamination
on $S'$ where $S_{j_1}$ is regarded as a double cover of $S'$ via $M$.
Then there exists a geometrically infinite Kleinian group $(\G, \psi)$ with the 
following properties on the boundary of $QH(G)$ 
in $AH(G)$.

\begin{enumerate}
\item The hyperbolic 3-manifold $\mathbf{H}^3/\G$ is topologically 
tame, i.e., there is a compact 3-manifold $M'$ whose interior is 
homeomorphic to $\mathbf{H}^3/\G$.
\item There is a homeomorphism $f : M \rightarrow M'$ inducing (a conjugate of)
the isomorphism $\psi : G \rightarrow \G$ when we identify $G$ 
with $\pi_1(M)$ and $\G$ with $\pi_1(M')$.

\item For each $i_k = i_1, \ldots , i_q$, the end of $\mathbf{H}^3/\G$ corresponding to $f(S_{i_k})$ 
is geometrically finite and its marked conformal structure at infinity
coincides with
$f_*(m_k)$.
\item For each $j_k = j_1, \ldots , j_p$, the lamination $f(\lambda_{j_k})$ 
is unrealisable in $M'$.
(See Remark \ref{rem}-(2) below.)
\item
If $G$ is not free, then $\Gamma$ has no parabolic elements, and the end facing $f(S_{j_k})$ for $j_k=j_1, \dots , j_p$ is geometrically infinite.

\end{enumerate}
\end{thm}

We should remark the following.

\begin{rem}
\label{rem}
\begin{enumerate}
\item In Ohshika \cite{OhE}, we gave a sufficient condition for a collection of 
simple closed curves to be made parabolic elements in a geometrically finite 
group on the boundary of quasi-conformal deformation space.
The theorem above can be regarded as a geometrically infinite 
version of this result.
In \cite{Ohpr}, which is a sequel to this one, a generalisation of our
main theorem
here allowing laminations to be disconnected and parabolic elements to exist is
given.

\item

In the theorem, we have only stated that $f(\lambda_{j_k})$ is unrealisable in $M'$.
To show that $f(\lambda_{j_k})$ really represents an ending lamination in
$M'$, we need to show the uniqueness of unrealisable laminations (in the Masur
domain) up
to diffeomorphism homotopic to the identity.
As was explained before, this fact is not easy to prove particularly in the case when the group is free, and its proof can be found in Ohshika \cite{Ohpr}.
An alternative proof of this fact in the case of free groups can be found in
Namazi-Souto \cite{NS}.

The tameness of the limit group $\Gamma$ needs the result of Brock-Souto \cite{BS}  and Brock-Bromberg-Evans-Souto \cite{BBES} or a general resolution of Marden's conjecture by Agol \cite{Ag} and Calegari-Gabai \cite{CG}.
In other cases, the unrealisability of the $f(\lambda_{j_k})$ implies the tameness.
%
\end{enumerate}
\end{rem}

%

We  present here an outline of the proof of the main theorem.
We shall construct a sequence of quasi-conformal deformations of $G$ using the Ahlfors-Bers
map so that for each boundary component on which a measured lamination is given, the
corresponding conformal structures converge to the projective class of that lamination in the
Thurston compactification of the Teichm\"{u}ller space.
The crucial step of the proof is to show that such a sequence converges after passing to a
subsequence.
By applying Bonahon's theory of
characteristic compression body to a compact core of
$\mathbf{H}^3/G$, we shall express $G$ using amalgamated free products and HNN extensions composed of function groups and
freely indecomposable groups.
Using Kleineidam-Souto's theorem, we see that if we restrict $\phi_i$ to a factor which is a
function group, then the sequence converges after taking a subsequence.
On the other hand, we can make use of Morgan-Shalen's interpretation of Thurston's theory by using the
language of $\R$-trees to show that the restriction of $\phi_i$ to a freely indecomposable factor also converges.
It will remain to show that these imply that the original groups, which can be obtained by
amalgamated free products and HNN extensions from these groups, also converges.
This will be shown using the fact the amalgamating subgroups converge to a quasi-Fuchsian
group.

The author would like to express his gratitude to the referee for his/her comments and suggestions which, the author hopes, have made this paper clearer and more readable.

\section{Preliminaries}
Kleinian groups are discrete subgroups of $PSL_2\C=SL_2 \C/\{\pm E\}$.
We regard Kleinian groups as acting on $\H^3$ by isometries and on its sphere at infinity
$S^2_\infty$, which is identified with the Riemann sphere, by conformal automorphisms.
Throughout this paper, we assume Kleinian groups to be
{\em finitely  generated and torsion free}.
A non-identity element of a (torsion-free) Kleinian group  is either  
loxodromic  or parabolic.
For a Kleinian group $G$, its limit set on $S^2_\infty$ is denoted by $\Lambda_G$.
The complement of $\Lambda_G$ in $S^2_\infty$ is called the domain of discontinuity 
of $G$, and is denoted by $\Omega_G$, on which  $G$ acts properly discontinuously.
The quotient $\Omega_G/G$ is a Riemann surface, which is known to be of 
finite type by work of Ahlfors.

For a hyperbolic 3-manifold $N$, its convex submanifold that is minimal among 
those which are deformation retracts of $N$ is called the {\em convex 
core} of $N$.
We say that a Kleinian group $G$ is {\em geometrically finite} when 
the convex core of $\H^3/G$ has finite volume.
A Kleinian group $G$ and $\H^3/G$ are said to be {\em topologically 
tame} when $\H^3/G$ is homeomorphic to the interior of a compact 
3-manifold.
As was shown by Marden \cite{Ma}, geometrically finite Kleinian groups 
are topologically tame.
In particular, when $G$ is geometrically finite and has no parabolic 
elements, the Kleinian manifold $(\H^3 \cup \Omega_G)/G$ gives a 
compactification of $\H^3/G$.

For an open 3-manifold $N$, its compact 3-submanifold $C$ is said to 
be a {\em compact core} of $N$ if the inclusion of $C$ into $N$ is a homotopy
equivalence.
Scott proved, in \cite{Sc}, that any irreducible open 3-manifold with finitely 
generated fundamental group has a compact core.
In particular, for a (finitely generated,
torsion-free) Kleinian group 
$G$, the quotient manifold $\H^3/G$ has a compact core.
It was proved by McCullough-Miller-Swarup \cite{MMS} that for any two compact cores $C_1, C_2$ of an
irreducible
$3$-manifold 
$M$, there is a homeomorphism from $C_1$ to $C_2$ inducing an inner-automorphism of $\pi_1(M) \cong
\pi_1(C_1) \cong \pi_1(C_2)$. 

A group is said to be {\em freely indecomposable} when it cannot be 
decomposed into a non-trivial free product, otherwise freely 
decomposable.
When a 3-manifold has fundamental group which is finitely generated and freely 
indecomposable, its compact core is boundary-irreducible, i.e., every 
boundary component is incompressible.

Let $N$ be a hyperbolic 3-manifold without cusp having finitely generated fundamental 
group, and $C$ its compact core.
Then the ends of $N$ correspond one-to-one to the
boundary components of $C$, since for each boundary component $S$ of $C$, 
there is only one end contained in the component of $N \setminus C$ 
touching $S$.
We say that the end as above {\em faces} $S$ in this situation.

An end of a hyperbolic 3-manifold $N$ without cusp is said to be geometrically 
finite when it has a neighbourhood intersecting no closed geodesics 
in $N$, otherwise  geometrically infinite.
The hyperbolic manifold $N$ is geometrically finite if and only if 
all the ends of $N$ are geometrically finite.
For each geometrically finite end of $N=\H^3/G$, there is a unique 
component $\Sigma$ of $\Omega_G/G$ which is attached to 
the end $e$ in 
the Kleinian manifold $(\H^3 \cup \Omega_G)/G$.
In this situation, the conformal structure of $\Sigma$ is called 
the conformal structure at infinity of the end $e$.

Let $S$ be a closed hyperbolic surface.
A {\em geodesic lamination} on $S$ is a  closed subset of $S$ consisting of 
disjoint simple geodesics.
The geodesics constituting a geodesic lamination are called the 
leaves.
A geodesic lamination endowed with a transverse measure on arcs which 
is invariant by  translations along leaves is called a 
{\em measured lamination}.
(We assume the measure to be non-zero for any arc with its interior intersecting the leaves.)
The set of measured laminations on $S$ endowed with a weak topology 
with respect to the measures on transverse arcs is denoted by 
$\mathcal{ML}(S)$ and called the measured lamination space of $S$.
The space which we obtain by taking the quotient of $\mathcal{ML}(S) \setminus
\{\emptyset\}$, 
identifying a measured
lamination with another obtained  by multiplying the transverse measure by a scalar, is called
the projective measured lamination  space and denoted by $\mathcal{PML}(S)$.
We call points in $\mathcal{PML}(S)$ projective laminations.
A measured lamination or a projective lamination is said to be 
{\em maximal} when it is not a proper sublamination of another measured 
lamination or projective lamination.
The underlying geodesic lamination of a measured lamination or a 
projective lamination is called its {\em support}.
A measured lamination or a projective lamination is said to be 
{\em connected} when its support is connected as a subset of $S$.

Here we shall give a remark on correspondence of our terminology with that  as 
was used in Otal
\cite{Ot} and Kleineidam-Souto
\cite{KS}.
There, they defined a geodesic lamination or a measured lamination to be {\em
arational} when its complementary regions are all simply connected.
{\em A measured lamination is arational if and only if it is maximal and
connected in our sense.}
They also use the term ``minimal lamination", which is equivalent to
``connected lamination" in the case of measured laminations.
Arational measured laminations are automatically minimal.

Thurston defined a natural compactification of the Teichm\"{u}ller space 
$\T(S)$ whose boundary is identified with $\PL(S)$.
(Recall that the Teichm\"{u}ller space is the space of marked 
conformal structures on $S$ modulo isotopic equivalence.
This can be regarded as the space of marked hyperbolic structures since
there is a one-to-one correspondence between the conformal structures and the hyperbolic
structures.)
 Let $\s$ be the set of isotopy classes of simple closed curves on
$S$. We consider the space of non-negative functions from $\s$ to $\R_+$, 
denoted by $\R_+^\s$, and its projectivisation $P\R_+^\s$.
The Teichm\"{u}ller space $\T(S)$ is embedded into $\R_+^\s$ by 
taking $g \in \T(S)$ to a function whose value at $s \in \s$ is the 
length of the closed geodesic representing $s$ with respect to $g$.
Likewise, the measured lamination space $\ML(S)$ is embedded in $\R_+^\s$ by 
taking $\lambda$ to a function whose value at $s$ is $i(\lambda, 
s)$, where $i(\lambda, s)$ is defined to be the infimum of 
$\int_\sigma d\mu$ for  the transverse measure $\mu$ of $\lambda$ and $\sigma$ 
ranging over all the simple closed curves representing $s$.
The Thurston compactification coincides with the compactification of 
the image of $\T(S)$ in $P\R_+^\s$ for these embeddings.
Refer to Fathi-Laudenbach-Po\'{e}naru \cite{FLP} for more details.

Let us now consider the situation where a hyperbolic surface $S$ is 
identified with a boundary component of a compact core of $C$ a hyperbolic 
3-manifold $N$.

\begin{dfn}[(Masur domain)]
We define a subset $\mathcal{C}(S)$ of the measured lamination space 
$\mathcal{ML}(S)$ to be the set of weighted disjoint unions of simple closed
curves bounding compressing discs in
$C$. Then the {\em Masur domain} $\M(S)$ is defined as
$$\M(S) = \{\lambda \in \mathcal{ML}(S) | i(\lambda, c) 
\neq 0 \text{ for all } c \in \overline{\mathcal{C}}(S)\},$$
where $\overline{\mathcal{C}}(S)$ denotes the closure of $\mathcal{C}(S)$ in $\mathcal{ML}(S)$.
Also, we call the image of the Masur domain in the projective 
lamination space the {\em projectivised Masur domain}.
\end{dfn}

We should note that the definition of the Masur domain in Otal 
\cite{Ot} differs from the one we gave above in the case when $C$ is 
a ``small compression body'' obtained as a boundary connected sum of two $I$-bundles over
closed surfaces. This difference does not matter in our setting, for we only deal
with maximal and connected measured laminations in the Masur domain.


Let $(S,m)$ be a closed hyperbolic surface and $N$ a hyperbolic 
3-manifold.
A continuous map $f: S \rightarrow N$ is said to be a pleated surface 
when
\begin{enumerate}
\item the length metric on $S$ induced from that on $N$ coincides 
with that of  $m$, and 
\item there is a geodesic lamination $\ell$ on $(S,m)$ such that both 
$f|\ell$ and $f|(S \setminus \ell)$ are totally geodesic.
\end{enumerate}
A geodesic lamination or a measured lamination on $S$ which  a 
pleated surface $f$ maps totally geodesically is said to be {\em realised} 
by $f$.

With regard to a measured lamination $\lambda$ on a surface $S$ contained in a hyperbolic $3$-manifold $N$, we say that $\lambda$ is {\em realisable} if there is a pleated surface homotopic to the inclusion of $S$ realising $\lambda$, and {\em unrealisable} otherwise.

An ending lamination for an end of a topologically tame 
3-manifold is defined as follows.
Let $N$ be a topologically tame hyperbolic 
3-manifold without cusps, having a compact core $C$ such that $N\setminus C
\cong \partial C \times \R$.
Let $e$ be a geometrically infinite end of $N$ facing a boundary component $S$ of $C$.
Then, as was shown by Canary \cite{CaJ}, there is a sequence of simple 
closed curves $\{c_j\}$ on $S$, whose projective classes converge to a 
projective lamination $[\lambda]$ contained in the projectivised Masur domain of 
$S$  regarded as 
the boundary component of $C$, such that there are closed geodesics $c_j^*$ 
homotopic to $c_j$ in $\overline{N \setminus C}$ tending to the end $e$ as $j
\rightarrow 
\infty$.
In this situation, the support of $\lambda$ is defined to be the ending
lamination of $e$, and 
the projective lamination
$[\lambda]$ or its  representative $\lambda$ is
said to represent the 
ending lamination of $e$.
It is known that the ending lamination of a topologically tame end is unique.
(This was proved by Bonahon in the case when $C$ is boundary-irreducible.
Canary showed that this can be generalised to the case when $C$ may be
boundary-reducible.)
An ending lamination is always unrealisable.
Conversely it has been proved that if  $N$ is an algebraic limit of geometrically finite groups, then
for a
maximal and connected measured lamination $\lambda$ in $\mathcal{M}(S)$  which is
unrealisable, there is a
homeomorphism $h$ of $C$ acting on $\pi_1(C)$ as an inner-automorphism such that 
$h(\lambda)$ represents the ending lamination of the end facing $S$.
(See Ohshika \cite{Ohpr} or Namazi-Souto \cite{NS}.)
We shall not use this fact in the proof of our main theorem.
Our result using this is given in \cite{Ohpr}.

A {\em compression body} is a connected compact 
3-manifold $V$ with $\partial V= \partial_e V \sqcup\partial_i V$, 
where we call $\partial_e V$ the exterior boundary, and $\partial_i V$ 
the interior boundary, 
such that no components of $\partial_i V$ are 
spheres  and  $V$ is obtained by attaching disjoint $1$-handles to a 
product neighbourhood of $\partial_i V$.
As an exceptional case, we also regard handlebodies as compression 
bodies whose interior boundaries are empty.
The exterior boundary $\partial_e V$ consists of only one component, which is compressible.
The interior boundary $\partial_i V$ is incompressible.
We do not allow the attached 1-handles to be empty, and accordingly, do 
not regard trivial $I$-bundles over closed surfaces as compression 
bodies.
For a disjoint union of compression bodies, we denote by $\partial_e$ and $\partial_i$ the
union of exterior boundaries and the union of interior boundaries
respectively.

Bonahon showed in \cite{BoE} that for any compact irreducible 
3-manifold $C$, there exists a disjoint union $V$ of  compression bodies in $C$, unique up 
to isotopy, such that the exterior boundary $\partial_e V $ is the 
union of all the compressible components of $\partial C$, the interior 
boundary $\partial_i V$ lies in $\I C$, and $\overline{C \setminus  V}$ is
irreducible 
and boundary-irreducible.
(When $C$ is boundary-irreducible, we set $V$ to be empty.)
Such a union of compression bodies is called the {\em characteristic compression body} 
of $C$.

For a Kleinian group $G$, its deformation space, denoted by $AH(G)$, 
is the set of equivalence classes of pairs $(\Gamma, \psi)$, where $\psi$ is a faithful 
discrete representation of $G$ into $PSL_2\C$ taking parabolic 
elements of $G$ to parabolic elements, and $\Gamma$ its image.
It should be noted that we allow loxodromic elements of $G$ to be mapped to parabolic elements
of $\Gamma$.
Two pairs are identified in $AH(G)$ when they are conjugate 
by an element of $PSL_2\C$.
We endow $AH(G)$ with the quotient topology induced from the space of 
representations of $G$ into $PGL_2\C$.

A Kleinian group $(\Gamma, \psi) \in AH(G)$ is called a 
quasi-conformal deformation of $G$ when there exists a quasi-conformal 
homeomorphism $f: S^2_\infty \rightarrow S^2_\infty$ such that 
$f\gamma f^{-1}= \psi(\gamma)$ for every $\gamma \in G$ as actions on 
$S^2_\infty$.
The subset of $AH(G)$ consisting of quasi-conformal deformations of 
$G$ is denoted by $QH(G)$.
It is known that $QH(G)$ is an open subset of $AH(G)$ when $G$ is geometrically 
finite.

Consider a point $m$ in the Teichm\"{u}ller space $\T(\Omega_G/G)$.
The point $m$ can be realised by a quasi-conformal homeomorphism preserving the markings from 
the surface $\Omega_G/G$ with the original marked conformal structure to 
$\Omega_G/G$ with the marked conformal structure $m$.
It is known, by work of Bers, that there exists a quasi-conformal homeomorphism 
on $S^2_\infty$ inducing a   quasi-conformal deformation of $G$,
whose restriction to $\Omega_G$ induces the quasi-conformal map above.
Moreover, such a quasi-conformal deformation is uniquely determined as an
element in $QH(G)$.
Therefore, there is a map $q : \T(\Omega_G/G) \rightarrow QH(G)$ 
taking $m$ to the corresponding quasi-conformal deformation.
This map $q$ is known to be a covering map, and is called the 
{\em Ahlfors-Bers map}.

For an isomorphism $\psi$ from a Kleinian group $G$ to another 
Kleinian group $\Gamma$, there is 
a homotopically unique homotopy equivalence from $\H^3/G$ to 
$\H^3/\Gamma$.
We denote such a homotopy equivalence by the same letter as the 
isomorphism but in the upper case, e.g., $\Psi$ for the isomorphism $\psi$.

An {\em $\R$-tree} is a complete metric space in which for any two points, there
is a unique simple arc connecting them, which is a geodesic.
When we say that a group $G$ acts on an $\R$-tree  $T$, we mean that for any $g
\in G$, there is an isometry $g: T \rightarrow T$ and the group operation
corresponds to the composition of isometries.
A group action on $T$ is said to have {\em small-edge stabilisers} when for every
non-trivial arc of $T$, its stabiliser is virtually abelian.

Let $N$ be a $3$-manifold and $T$ an $\R$-tree.
(We usually consider the case when $N$ is a universal cover of a compact $3$-manifold.)
A continuous map $f: N \rightarrow T$ is said to be {\em weakly transverse} when every $x
\in M$, there is a neighbourhood $U \cong D^2 \times I$ of $x$ such that $f|D^2
\times \{t\}$ is constant for every $t \in I$ and $f|\{pt.\} \times I$ is 
monotone in the weak sense, i.e., if $t$ lies in $[t_1,t_2] \subset I$, then $f(D^2 \times \{t\})$ lies in the (closed) segment between $f(D^2 \times \{t_1\})$ and $f(D^2 \times \{t_2\})$.
For a transverse map $f: N \rightarrow T$, we can consider a codimension-1
measured lamination on $N$ 
 whose support consists of points  where $f$ is not
locally constant, and whose local flow box is exactly $D^2 \times
I$-neighbourhood as above with the transverse measure equal to the metric induced on
$I$ from $T$ by $f$.
This lamination is called the the dual lamination of $f$ (or the lamination
dual to $f$).

An arc $\alpha$ in $N$ is said to be {\em monotonically transverse} to a lamination $L$ in $N$ when it is transverse to leaves and there is no subarc $\alpha': [a,b] \rightarrow N$ of $\alpha$ whose endpoints lie on the same leaf of $L$ and $\alpha'(a,b)$ lies in the complement of $L$, i.e, $\alpha$ cannot turn around in the complement and comes back to the same leaf.
We say that a weakly transverse map $f: N \rightarrow T$ is {\em strongly transverse} to $T$ if for any arc $\alpha$ monotonically transverse to the dual lamination, $f \circ \alpha$ is monotone in the weak sense, i.e., for any $t \in [t_1, t_2]$ the point $f \circ \alpha(t)$ lies on the segment joining $f \circ \alpha(t_1)$ and $f \circ \alpha(t_2)$.

As is obvious from the definition the difference between weakly transverse maps and strongly transverse maps lies in the fact that a weakly transverse map may have {\em folds}, i.e., there may be a monotonically transverse arc in $\tilde{C}$ which is mapped to a non-monotone arc in the $\R$-tree.

Now consider the situation when $\tilde{C}$ is the universal covering of  a
$3$-manifold $C$, and let $G$ be the covering translation group.
Suppose that $G$ acts on an $\R$-tree $T$ isometrically.
Considering a handle decomposition of $C$, we can construct an equivariant
weakly transverse map from $\tilde{C}$ to $T$.
(Refer to Morgan-Shalen \cite{MS3}.)
Then the lamination dual to $f$ projects to a codimension-1 measured lamination
in $C$, which we also call the {\em lamination dual to} $f$.

From the lamination dual to a weakly transverse $f$ above, we can construct an $\mathbf{R}$-tree $T'$ by taking the completion of the leaf space, and $\pi_1(C)$ acts on $T'$ by isometries.
Nevertheless, this tree $T'$ {\em does not} necessarily coincide with the original tree $T$.
If $f$ is {\em strongly} transverse, there is an isometric embedding from $T'$ to $T$ which is  equivariant under the action of $\pi_1(C)$.
Therefore, if $T$ is assumed to be minimal, then the two $\R$-trees coincide in this case.

\section{Some results due to Thurston and Morgan-Shalen}
In this section, we shall present three theorems all of which are 
essentially due to Thurston.
%
The first of them is the following.

\begin{thm}
\label{Thurston}
Let $G$ be a freely indecomposable Kleinian group without parabolic 
elements.
Let $C$ be a compact core of $\mathbf{H}^3/G$.
Suppose that we have a sequence $\{(\Gamma_i, \psi_i) \in AH(G)\}$ 
which does not have a convergent subsequence.

Then there are a sequence of disjoint unions of  essential annuli 
$\{A_i^1\sqcup \ldots \sqcup A_i^k\}$ 
properly embedded in $C$, whose number $k$ is constant with respect to 
$i$,  and positive weights $w_i^1, \ldots , w_i^k$, satisfying the following.
The weighted union of annuli $w_i^1 A_i^1 \sqcup \dots \sqcup w_i^k
A_i^k$ converges to a codimension-1 measured lamination in $C$. 
Let $\{r_ic_i\}, \{s_i d_i\}$ be  sequences of weighted simple closed curves on $\partial C$
converging to measured laminations in
$\mathcal{ML}(\partial C)$.
Then,
$$\lim_{i \rightarrow \infty}
\frac{\mathrm{length}(\Psi_i(r_ic_i))}{\mathrm{length}(\Psi_i(s_i d_i))}=\lim_{i
\rightarrow
\infty}\frac{ w_i^1
i(A^1_i,r_i c_i )+
\dots 
+w_i^ki(A^k_i, r_i c_i)}{w_i^1
i(A^1_i,s_i d_i )+
\dots 
+w_i^ki(A^k_i, s_i d_i)}$$ after taking a subsequence of
$\{(\Gamma_i,\psi_i)\}$ whenever the limit of the right hand does not have
form of $0/0$.
Here for a weighted simple closed curve $r_ic_i$, the length denoted by $\mathrm{length}(\Psi_i(r_i c_i))$ is defined to be $r_i\mathrm{length}(c_i^*)$ for the closed geodesic $c_i^*$ homotopic to $\Psi_i(c_i)$ in $\mathbf{H}^3/\Gamma_i$.
Moreover, if
$\lim_{i
\rightarrow
\infty} w_i^1 i(A^1_i,r_i c_i )+
\dots 
+w_i^ki(A^k_i, r_i c_i) >0$, then $\lim_{i \rightarrow \infty}
\mathrm{length}(\Psi_i(r_i c_i)) = \infty$.
\end{thm}

\begin{rem}
This theorem {\em does not} say that $C$ contains a measured lamination geometrically realising the rescaled limit of  $(\Gamma_i, \phi_i)$, which is an action of $G$ on an $\R$-tree.
(The limit in this sense is the same as that in $P\R^\mathcal{S}$, where $\mathcal{S}$ denotes the set of free homotopy classes of essential closed curves in $C$ and we regard $(\Gamma_i, \psi_i)$ as a point there by letting an $s$-coordinate be the translation length of $\psi_i(s)$ for every $s \in \mathcal{S}$.)
The point is that in our conclusion, we only consider weighted simple closed curves on $\partial C$.
We do not claim that the same kind of equality holds for every homotopy class of  essential closed curves in $C$.
In the corresponding theorem stated without proof (Theorem 3.1) in our previous paper \cite{OhZ}, this distinction was not made clear, and accordingly the statement was misleading.
What was really needed there was nothing but Theorem \ref{Thurston} above.

Actually, there is an example of $\{(\Gamma_i, \psi_i)\}$ diverging in $AH(G)$ whose rescaled limit cannot be realised by a measured lamination in $C$.
(The author was informed that such a phenomenon was first observed by J-P. Otal.)
This can be constructed as follows.
First assume that the characteristic pair of $C$ contains a solid torus component $V$ such that $\partial V \setminus \partial C$ consists of more than three open annuli whose core curves are homotopic in $V$ to the axis of $V$.
For simplicity, we assume that $\partial V \setminus \partial C$ consists of four open annuli, whose closures we denote by $A_1, A_2, A_3,$ and $A_4$ assuming they lie in this order on $\partial V$.
We can construct a $3$-manifold $C'$ homotopy equivalent to $M$ by pasting $C \setminus V$ to $V$ by interchanging what was pasted to $A_1$ to what was pasted to $A_2$.
Then the natural homotopy equivalence obtained by extending the identity map on $C \setminus V$ is not homotopic to a homeomorphism unless there is a symmetry of $C\setminus V$ interchanging $A_1$ with $A_2$.

Now, we consider a measured lamination in $C'$ consisting of only one essential annulus $A$ with the unit weight  properly embedded in $(V, V \cap \partial C')$ separating $A_2, A_3$ from $A_4, A_1$.
There is an action of $\pi_1(C)$ on a simplicial tree $T$ corresponding to the preimage of $B$ in the universal cover of $C'$.
It turns out that this action cannot be realised by any measured lamination in $C$.
In fact, any essential annulus in $V$ separates in $C$ either what is pasted to $A_21$ (under the rule of pasting for $C$)  from what is  pasted to $A_3$ or what is pasted to $A_2$ from what is pasted to $A_4$ whereas $A$ separates neither of them (under the rule of pasting for $C'$) in $C'$.

On the other hand, this annulus $A$ with unit weight appears obviously as a rescaled limit of geometrically finite  hyperbolic structure without cusps in $\mathrm{Int} C'$.
Since $C'$ is obtained from $C$ by what Canary-McCullough \cite{CM} called a primitive shuffling, the subspaces, $CC(C)$ in $AH(C)$ consisting of geometrically finite hyperbolic structures without cusps in $\mathrm{Int}C$ and $CC(C')$ consisting of those in $\mathrm{Int}C'$, bump each other, and  we can easily see that this rescaled limit is realised as that of groups in $CC(C)$ as follows.
A sequence of $CC(C)$ whose rescaled Gromov limit is $A$ with unit weight is obtained by pinching the conformal structure at infinity along two simple closed curves corresponding to the components of $\partial A$ on $\partial C'$ at the same speed.
This sequence can be approximated by bumping points of $CC(C)$ and $CC(C')$ corresponding to a geometrically finite structure with a $\Z$-cusp.
Thus, this gives a rescaled limit which cannot be realised by a measured lamination in $C$.

We should note that in this example a sequence of weighted annuli having the properties of our theorem can be obtained as a disjoint union of two essential annuli $A' \sqcup A''$ in $V \subset C$ such that $A'$ separates $A_1$ from the rest of annuli whereas $A''$ separates $A_3$ from the rest, which is constant with respect to $i$.
\end{rem}

\begin{proof}[Proof of Theorem \ref{Thurston}]
Although this theorem is essentially due to Thurston,
we shall give its proof here since a detailed proof can be found in no
references.
Here we adopt an argument using $\R$-trees based on the results of
Morgan-Shalen \cite{MS3} and Skora \cite{Sk} combined with the efficiency of pleated surfaces proved by Thurston.
There may be an alternative approach using the results of Thurston in \cite{Th2} and \cite{Th3}.

Since we assumed that $G$ is freely indecomposable, 
the compact core $C$ of $\mathbf{H}^3/G$ is boundary-irreducible.
 Consider the characteristic pair $W$ of $C$.
 (See Jaco-Shalen \cite{JS} and Johannson \cite{Jo} for the 
definition of characteristic pair.)
Since $C$ is atoroidal and has no torus boundary components, each component of
its characteristic pair is either an $I$-pair or a solid torus.
The characteristic pair has a property that every  essential annulus can be 
homotoped (as a pair of maps) into it. 
If a frontier component of an $I$-pair in the characteristic pair is
homotopic
to a frontier component of another component which is a solid tours, then an
annulus homotopic to such a frontier component has two ways to be carried by the
characteristic pair.
In the following argument, we always regard such an annulus as being carried
by a solid torus component.

We  use the following reformulation of Thurston's work by Morgan-Shalen
\cite{MS3} using
$\R$-trees.

\begin{ththm}[Morgan-Shalen]
\label{Morgan-Shalen}
Suppose that $\{(\Gamma_i, \psi_i)\}$ does not converge in $AH(G)$ even after 
passing to a subsequence.
Then there are an action of $G$ on an $\R$-tree $T$ with small edge-stabilisers
and
a weakly transverse $G$-equivariant map $f$ from the universal cover $\tilde{C}$
of $C$ to $T$ with the following conditions.
\begin{enumerate}
\item There is a sequence of positive real numbers $\{\epsilon_i\}$ going to $0$ such that
\linebreak
$\lim_{i \rightarrow \infty}\epsilon_i\l_{\H^3/\Gamma_i}(\psi_i(\gamma))= \l_T(\gamma)$ for every
$\gamma
\in G$, where $\l_T$ denotes the translation length on $T$.
(We say that $\{(\Gamma_i, \psi_i)\}$ diverges to this $\R$-tree action in this situation.)
\item The dual lamination $L$ on $C$ of $f$ is carried by an incompressible branched surface $B$.
\item The branched surface $B$, hence also the lamination $L$, can be properly
isotoped into the characteristic pair of $C$.
\item 
The branched surface $B$ carries only surfaces with Euler number $0$.
\end{enumerate}
\end{ththm}


To derive Theorem \ref{Thurston} from the theorem above, we need to deform the weakly transverse map $f$ so that if it is restricted to each boundary component $S$ of $C$ it becomes a strongly transverse map to a subtree of $T$ invariant under $\pi_1(S)$.
For that, first we need to refine Theorem \ref{Morgan-Shalen} as follows to make it suitable
for our purpose.

\begin{cor}
\label{refined Morgan-Shalen}
Let $S$ be a component of
$\partial C$, and $\tilde{C}$ the universal cover of $C$.
In the situation of Theorem \ref{Morgan-Shalen}, let
$T_S$ be the minimal invariant sub-$\R$-tree  of $T$ with respect to
$\pi_1(S)$ regarded as a subgroup of $G$.
Let $\tilde{S}$ be a component of the preimage of $S$ in $\tilde{C}$ which is
invariant under $\pi_1(S) \subset G$ acting as covering translations.
Then, we can choose $f$ in the conclusion of  Theorem \ref{Morgan-Shalen} so that
$f|\tilde{S}$ is a weakly transverse map into $T_S$.
\end{cor}
\begin{proof}
We first review an outline of the construction of $f$ in the proof of  Theorem
\ref{Morgan-Shalen} provided that an action of $G$ on an $\R$-tree $T$ as in the theorem was given.
(The existence of such an action was proved in Morgan-Shalen \cite{MS}.
Also, there are alternative proofs by Paulin
\cite{Pa} and Bestvina \cite{Be}.)
First fix a handle decomposition of $C$.
We can find a weakly transverse map $g: \tilde{C} \rightarrow T$ whose dual lamination $\Lambda$ on $C$ is normal with respect to the handle decomposition.
Let $B$ be a branched surface carrying $\Lambda$.
We make $B$ incompressible by performing surgery which corresponds to
moving $g$ by a homotopy, preserving the condition that $B$ carries the dual
lamination of $g$. 
In this process, we homotope $g$ in neighbourhoods of  compressing discs  and
discs of contact to
compress $B$, and in a ball bounded by an inessential sphere carried by $B$ to
remove a spherical part of the boundary.
Similarly we can make $B$ boundary-incompressible by the boundary-compression
(along either a boundary-compression disc or a semi-disc of contact) and 
the removal of inessential part carrying a boundary-parallel disc.
(Although this process needs a limit argument, the essence of the
operation consists of 
these two.)
Let $f$ be the thus modified weakly transverse equivariant map.
It turns out that the assumption that the action on $T$ has small edge-stabilisers implies that the branched surface $B$ carrying the dual lamination can carry only surfaces with null Euler characteristics.
By the Jaco-Shalen-Johannson theory, it follows that the dual lamination can be
isotoped into the characteristic pair then.

Now we shall see how we can make $f$ satisfy the condition of our corollary.
In the construction of the weakly transverse map $g$, it is clear that we can make $g(\tilde{S})$ lie in $T_S$.
In the process of the modification of the map, the only part that affects the
restriction to $\tilde{S}$ is that of boundary-compression, removing
semi-discs of contact and removing parts carrying a boundary-parallel discs.
We can choose a boundary-compressing semi-disc $\Delta$ so that its intersection with the dual lamination consists of arcs parallel to the arc $\overline{\Delta \setminus \partial C}$.
Then, since in the lift $\tilde{\Delta}$ of $\Delta$ attached to $\tilde{S}$,
every point has a point on $\tilde{S}$ with the
same image by $g$ in $T$, we see that $\tilde{\Delta}$ is mapped into $T_S$.
Therefore, after the boundary-compression, $\tilde{S}$ remains to be mapped
into $T_S$.
If the limit is taken within the maps with this property, the limit map has
the same property. 
The same argument applies to semi-discs of contact which intersect the dual
lamination.
(If it is disjoint from the dual lamination, we can remove it without changing
the weakly transverse map.)
It is easy to see that removing a part carrying a boundary-parallel disc does not
change the condition that $\tilde{S}$ is mapped into $T_S$ since a lift $\tilde{D}$ of such a disc can be taken to
be mapped into $T_S$ and the resulting map takes the disc on $\tilde{S}$ parallel to $\tilde{D}$ to the image of  $\tilde{D}$. 
\end{proof}

\begin{lem}
\label{Skora}
We can homotope $f$ equivariantly only in neighbourhoods of the preimages of the compact leaves of $L|\partial C$ so that for every component $S$ of $\partial C$ with $L \cap S \neq \emptyset$, the restriction
$f|\tilde{S}$ is a strongly transverse map onto $T_S$, hence induces an isomorphism, i.e., an  isometry from   the $\R$-tree obtained as the completion of the leaf space of the measured lamination
dual to $f|\tilde{S}$, to  $T_S$ which is equivariant under the action of $\pi_1(S)$.
(If $S$ is disjoint from $L$, then $T_S$ consists of only one point and $f$ maps $\tilde{S}$ to the point.) 
This process corresponds to identifying  or
removing some of the compact leaves of $L|\partial C$.
For the thus obtained transverse map $f$, the condition of Theorem \ref{Morgan-Shalen} remains valid.
\end{lem}

\begin{proof}
Let $S$ be a component of $\partial C$ with $L \cap S \neq \emptyset$.
For a weakly transverse map $f$, we denote the lamination on $S$ dual to $f|\tilde{S}$
by $\lambda_f$ and that on $\tilde{S}$ by $\tilde{\lambda}_f$.
By Proposition 3.1 in Skora  \cite{Sk}, if $f|\tilde{S}$ does not fold at an edge point, then  $f|\tilde{S}$ is strongly transverse to $T_S$, hence induces an isomorphism between the action of $\pi_1(S)$ on the $\R$-tree obtained from the completion of the leaf space of the lamination $\lambda_f$ and the one on $T_S$ which is the restriction of the given action on $T$.
Furthermore, if an edge point $x$ corresponds to (a lift of) a leaf of minimal component which is not  a closed leaf, then $f|\tilde{S}$ cannot fold at $x$ since the leaf containing $x$ is dense in the minimal component.
(This means that if $y$ does not lie on a boundary leaf, for a short arc $I$
transverse to the lamination
at $y$, the map $f|I$ is monotone.)

Suppose now that $x$ corresponds to (a lift of) a compact leaf of $\lambda_f$.
Note that the dual lamination $\lambda_f$ may have parallel compact leaves
which correspond to one compact leaf as a geodesic lamination.
Consider all the leaves of $\lambda_f$ that are parallel to the one corresponding to $x$.
Let $A$ be the smallest annulus on $S$ containing all these leaves.
Then both components of $\partial A$ are also leaves.
Let $\tilde{A}$ be a lift of $A$ in $\tilde{S}$.
Take points $a,b$ on $\partial \tilde{A}$, one on each boundary component.
Let $d$ be $d_T(f(a), f(b))$.
If $d >0$, then we can homotope $f$ equivariantly only within the translates of
$\tilde{A}$ in such a way that  $f|\tilde{A}$ becomes
monotone in the transverse direction and the total transverse measure of the
leaves in $\tilde{A}$ induced by $f$ is equal to $d$. 
If $d=0$, we can homotope $f$ similarly within the translates of $\tilde{A}$
and remove all the leaves in the translates of $\tilde{A}$ from the dual
lamination.
By repeating this operation for every compact leaf, we can make $f$ not
fold at an edge point.
Thus we have shown, by Skora's Proposition 3.1, we can homotope $f$ so that $f|\tilde{S}$
is a strongly transverse map to $T_S$.
The map is surjective since we assumed $T_S$ to be minimal.

In each step of this operation, the codimension-one measured lamination of $C$
dual to this new $f$ may change.
By this operation, some leaves which are annuli contained in a 
component $V$ of the characteristic pair are
glued to other annuli in another component $V'$ and are pushed into the
interior to form
new annular leaves.
If these annuli are essential, then they are carried by the characteristic
pair, and
the conditions in Theorem
\ref{Morgan-Shalen} are satisfied.
If they are inessential, i.e., boundary-parallel, then we can further homotope
$f$ and remove these leaves as before.
We perform this modification of $f$ for each component of $\partial C$.
Thus we have shown that the dual lamination of $f$ satisfies the conditions in
Theorem \ref{Morgan-Shalen}.
\end{proof}

Having proved the lemma above, we can prove Theorem \ref{Thurston} as follows.
Let $\{r_i c_i\}, \{s_i d_i\}$  be sequences of weighted simple closed curves
on $
\partial C$ as were given in the statement of the theorem.
Let $S_1, S_2$ be components on which $\{r_i c_i\}$ and $\{s_i d_i\}$ lie
respectively,
and 
$\lambda$ and $\mu$  measured laminations to which $\{r_i c_i\}$ and
$\{s_i d_i\}$ converge respectively.

At this point, we need another result by Thurston.
Let $S$ be a component of $\partial C$.
Then,  Theorem 3.3 of Thurston \cite{Th2} shows that
there are a continuous function $K: \mathcal{ML}(S) \rightarrow \mathbf{R}$ independent of
$i$ and a pleated surface $f_i : (S, m^S_i)  \rightarrow 
\mathbf{H}^3/\G_i$ homotopic to $\Psi_i|S$ for each $i$, such that for any $\lambda
\in 
\mathcal{ML}(S)$,
\begin{equation}
\label{1}
\mathrm{length}_{\mathbf{H}^3/\G_i}(\Psi_i(\lambda)) \leq 
\mathrm{length}_{m_i^S}(\lambda) \leq 
\mathrm{length}_{\mathbf{H}^3/\G_i}(\Psi_i(\lambda))  + K(\lambda).
\end{equation}
 This inequality is called the {\em efficiency of pleated surfaces}.
 It should be noted that this efficiency necessitates the assumption that the
surface is incompressible, which is true in our case because $C$ is
boundary-irreducible.
In the inequality, $\l_{\H^3/\G_i}(\Psi_i(\lambda))$ means, when $\lambda$ is connected, the
length of a realisation of
$\lambda$ by a pleated surface 
homotopic to $\Psi_i$ if it is realisable by such a pleated surface, and is 
set to be $0$ 
if it is not.
When $\lambda$ is disconnected, the length is defined to be the sum of the lengths of its
components.

Now, suppose first that both $S_1$ and $S_2$ intersect $L$.
Lemma \ref{Skora} implies that $\{m_i^{S_j}\}$ converges to a projective
lamination represented by $L \cap S_j$ for $j=1,2$ in the Thurston compactification
of the Teichm\"{u}ller space.
(See Skora \cite{Sk} and Otal \cite{OtA}.)
Take simple closed curves $\gamma$ on $S_1$ and $\delta$ on $S_2$ both intersecting $L$ essentially.
As was shown in Expos\'{e}  8 of Fathi-Laudenbach-Po\'{e}naru \cite{FLP}, this
implies that $$\lim_{i \rightarrow
\infty}\frac{\l_{m_i^{S_1}}(r_i c_i)}{\l_{m_i^{S_1}}(\gamma)}=\lim_{i \rightarrow
\infty}\frac{i(L\cap S_1, r_i c_i)}{i(L\cap S_1, \gamma)},$$ 
and
$$\lim_{i \rightarrow
\infty}\frac{\l_{m_i^{S_2}}(s_i d_i)}{\l_{m_i^{S_2}}(\delta)}=\lim_{i \rightarrow
\infty}\frac{i(L\cap S_2, s_i d_i)}{i(L\cap S_2, \delta)}.$$

Using the efficiency of pleated surfaces (\ref{1}), we have $$\lim_{i\rightarrow
\infty}\frac{\l_{\H^3/\Gamma_i}(\Psi_i(r_i c_i))}{\l_{\H^3/\Gamma_i}(\Psi_i(\gamma))}=\lim_{i
\rightarrow
\infty}\frac{\l_{m_i^{S_1}}(r_i c_i)}{\l_{m_i^{S_1}}(\gamma)},$$ and $$\lim_{i\rightarrow
\infty}\frac{\l_{\H^3/\Gamma_i}(\Psi_i(s_i d_i))}{\l_{\H^3/\Gamma_i}(\Psi_i(\delta))}=\lim_{i \rightarrow
\infty}\frac{\l_{m_i^{S_2}}(s_i d_i)}{\l_{m_i^{S_2}}(\delta)}.$$

On the other hand, by Theorem \ref{Morgan-Shalen},  there
is a sequence of positive numbers $\{\epsilon_i\}$ going to $0$  such that for every $g \in G$ , we have that $\epsilon_i \mathrm{length}_{\H^3/\Gamma_i}(\psi_i(g))$ converges to
the translation length of $g$ in $T$.
Therefore, for the simple closed curve $\gamma$ as above, we have $\lim_{i \rightarrow \infty}\epsilon_i \l_{\H^3/\Gamma_i}(\Psi_i(\gamma))=\l_T(\gamma)$
and by \linebreak Lemma
\ref{Skora}, this is equal to $i(\gamma, L \cap S_1)$ since $\gamma$ lies on $\partial C$.
Similarly, \linebreak$\lim_{i \rightarrow \infty} \epsilon_i
\l_{\H^3/\Gamma_i}(\psi_i(\delta))=i(\delta, L \cap S_2)$.
Combining these with the equations above, we see that  $\lim_{i
\rightarrow \infty}\epsilon_i\l_{\H^3/\Gamma_i}(\Psi_i(r_i c_i))= \lim_{i \rightarrow \infty}i(r_i c_i, L)$
and \linebreak
$\lim_{i \rightarrow \infty}\epsilon_i\l_{\H^3/\Gamma_i}(\Psi_i(s_i d_i))  =\lim_{i\rightarrow
\infty}i(s_i d_i, L \cap S_2)$.
Thus we have proved that 
$$\lim_{i \rightarrow \infty}\frac{\l_{\H^3/\Gamma_i}(\Psi_i(r_i c_i))}{\l_{\H^3/\Gamma_i}(\Psi_i(s_i d_i))}=
\lim_{i\rightarrow \infty}\frac{i(r_i c_i, L)}{i(s_i d_i,L)}.$$

Since $L$ is a codimension-1 measured lamination which is a limit of a
disjoint union of weighted essential annuli by Theorem \ref{Morgan-Shalen}, we
get the equality in the conclusion of our theorem under the assumption that both
$S_1$ and $S_2$ intersect $L$.
Since $\epsilon_i \rightarrow 0$ and $\lim_{i \rightarrow \infty} \epsilon_i \l_{\H^3/\Gamma_i}(\Psi_i(r_ic_i)) =
\lim_{i \rightarrow \infty}i(r_i c_i, L)$,
we also obtain the last sentence of the conclusion.
(This is true even when $S_2$ is disjoint from $L$.)

Next suppose that one of $S_1, S_2$, say $S_1$, is disjoint from $L$.
(We do not need to consider the case when both of $S_1, S_2$ are disjoint from
$L$, for the right hand limit has the form $0/0$ then.)
Then, each component of the preimage of  $S_1$ in $\tilde{C}$ is mapped to a
point by $f$.
Let $\tau$ be a train track carrying both $\lambda$ and the $c_i$ for large
$i$.
Then every component of the preimage of $\tau$ is mapped to a point by $f$.
By regarding the $c_i$ as lying on $\tau$ as a $C^1$-curve, and pulling back the image
of its preimage to $\H^3/\Gamma_i$ rescaled by $\epsilon_i$,
we have $\epsilon_i \l_{\H^3/\Gamma_i}(\psi_i(r_i c_i)) \rightarrow 0$.
(Compare with the argument in the chapitre 3 of Otal \cite{OtA}.)
Since we are not considering the case when the limit in the right hand side has
the form $0/0$, we can assume that $\lim_{i \rightarrow \infty} i(s_i d_i, L)>0$.
Then by the argument above, $\lim_{i \rightarrow \infty}\epsilon_i
\l_{\H^3/\Gamma_i}(\Psi_i(s_i d_i))
=\lim_{i \rightarrow \infty} i(s_i d_i, L)>0$.
Thus we have shown that both of the limits
$$\lim_{i \rightarrow \infty}\frac{\l_{\H^3/\Gamma_i}(\Psi_i(r_i c_i))}{\mathrm{length}_{\H^3/\Gamma_i}(\Psi_i(s_i d_i))}
\text{ and }
\lim_{i \rightarrow \infty}
\frac{i(r_i c_i, L)}{i(s_i d_i, L)}$$
are
equal to $0$ and that the equality holds.
\end{proof}

The following theorem appeared in Thurston \cite{Th2} as Theorem 2.2, 
whose proof can be found in the argument of Thurston \cite{Thm}.
\begin{thm}[Thurston]
\label{length to 0}
Let $S$ be a closed surface of genus greater than $1$.
Let $\{m_i\}$ be a sequence in the Teichm\"{u}ller space 
$\mathcal{T}(S)$, which converges to a projective lamination 
$[\lambda]$ (represented by a measured lamination $\lambda$) in the Thurston
compactification
of $\T(S)$. Then there is a sequence of weighted simple closed curves $\{w_i 
\gamma_i\}$ converging to the measured lamination $\lambda $ such that 
$\mathrm{length}_{m_i}(w_i\gamma_i) \rightarrow 0$ as $i \rightarrow 
\infty$.
\end{thm}

For a Kleinian group $G$, and a boundary component $S$ of a compact 
core $C$ of $\H^3/G$, we define the length function $\ell : \M(S) \times 
AH(G) \rightarrow \R$ as follows.
\begin{enumerate}
\item For $(\Gamma, \psi) \in AH(G)$, when $\lambda$ can be realised by a pleated surface
homotopic 
to $\Psi|S$, we set $\ell(\lambda, (\Gamma, \psi))$ to be 
the length of $\lambda$ on the pleated surface.
(It is known that the length does not depend on the choice of pleated surface realising
$\lambda$.)
\item When $\lambda$ is connected and there is no pleated surface homotopic to
$\Psi|S$
 realising $\lambda$, we set $\ell(\lambda, (\Gamma, \psi))$ to be $0$.
 \item When $\lambda$ is disconnected and is not realised, then we define
$\ell(\lambda, (\Gamma, \psi))$
to be the sum of the $\ell(\lambda_j, (\Gamma, \psi))$ for the components
$\lambda_j$ of $\lambda$.
\end{enumerate}

The  continuity of $\ell$ was first announced in Thurston
\cite{Th2} in the
 case when $S$ is incompressible.
In the subset consisting of maximal and connected laminations in
$\mathcal{ML}(S)$ (for incompressible $S$) its
proof  was given in Ohshika
\cite{OhG} based on Bonahon's work in \cite{BoA}.
  Brock
\cite{Br} gave a proof for the continuity of $\ell$ in the entire set
$\mathcal{ML}(S)
\times AH(G)$ still in the case when $S$ is incompressible.
The following theorem is a weak form of the continuity of $\ell$ in the
general case where $S$ may be compressible.
We can prove this using Bonahon's technique generalised in \cite{OhM}.
We shall sketch a proof in the appendix.

\begin{thm}
\label{continuity}
Let $\{(G_i,\phi_i)\} \in AH(G)$ be a sequence converging to $(\Gamma, \psi)$, 
 and  $\lambda_i$ measured laminations in $\mathcal{M}(S)$ converging to a
maximal and connected measured lamination $\mu
\in
\mathcal{M}(S)$. Suppose that $\ell(\lambda_i, (G_i, \phi_i))$ converges to $0$ as $i
\rightarrow \infty$. Then there is no pleated surface homotopic to $\Psi|S$ realising $\mu$ in
$\H^3/\Gamma$.
\end{thm}

\section{Convergence theorem}
The key step of the proof of Theorem \ref{main} is to give a 
sufficient condition for a sequence of quasi-conformal deformations 
to converge in the deformation space $AH(G)$, which is formulated 
as follows.

\begin{thm}
\label{converge}
Let $G$ be a geometrically finite Kleinian group without parabolic 
elements, which is not a quasi-Fuchsian group.
Let $M$ be a compact 3-manifold (Kleinian manifold) whose interior is identified
homeomorphically with 
$\mathbf{H}^3/G$, and $S_1, \ldots , S_m$ its boundary components.
Let $\lambda_{j_1}, \ldots , \lambda_{j_p}$ be maximal and connected measured laminations 
on boundary components $S_{j_1}, \ldots , S_{j_p}$ among $S_1, 
\ldots, S_m$, which are contained in the Masur domain.
 (We regard the Masur domain for an
incompressible surface as the entire measured lamination space.)
Consider a sequence of quasi-conformal deformations $\{(G_i, \phi_i)\}$ 
of $G$ induced from conformal structures at infinity  $g_1^i, 
\ldots , g_m^i$ on $S_1, \ldots , S_m$ such that  $\{g^i_{j_n}\}$ 
converges as $i \rightarrow \infty$ in the 
Thurston compactification of the Teichm\"{u}ller space to the 
projective lamination
$[\lambda_{j_n}]$ for $n=1,\ldots, p$, and for the rest of 
$k=1,\ldots , m$, the sequence $\{g^i_k\}$ stays in a compact set of 
the Teichm\"{u}ller space.
Then $\{(G_i,\phi_i)\}$ converges in $AH(G)$ after passing to a 
subsequence.
\end{thm}

In the special case when $G$ is a function group, this theorem was 
proved by Kleineidam and Souto in \cite{KS}.

\begin{thm}[Kleineidam-Souto]
\label{KS}
Let $G$ be a geometrically finite  Kleinian group without parabolic 
elements such that $\mathbf{H}^3/G$ is homeomorphic to the interior of 
a compression body $V$.
Let $S$ be the exterior boundary of $V$, and $\{m_i\}$ a sequence in 
the Teichm\"{u}ller space $\mathcal{T}(S)$ which converges in the 
Thurston compactification of $\T(S)$ to a maximal and connected 
projective lamination contained in the projectivised Masur domain.
Let $\{n_i\}$ be any sequence in the Teichm\"{u}ller space 
$\T(\partial_i V)$.
Let $(G_i, \phi_i)$ be a quasi-conformal deformation of $G$ which is 
given as the image by the Ahlfors-Bers map of $(m_i, n_i) \in 
\T(S) \times \T(\partial_i V)$.
Then $\{(G_i, \phi_i)\}$ converges in $AH(G)$ after passing to a 
subsequence.
\end{thm}

On the other hand, in \cite{OhI}, we proved a convergence theorem 
similar to Theorem \ref{converge} in a more general 
setting for freely indecomposable Kleinian groups, which
implies, in particular, Theorem \ref{converge} in the case 
when $G$ is freely indecomposable.
(This was based on Thurston's theorem which is the same as
Theorem
\ref{Thurston} in the present paper.)
Note that the assumption in Theorem \ref{converge} that measured laminations are contained in 
the Masur domains is void in this case.
This theorem in \cite{OhI} also covers the case when $G$ is a quasi-Fuchsian 
group with the following extra assumption:
In the case when $G$ is a quasi-Fuchsian group corresponding to an orientable
surface and $p=2$, we further
 assume
 that   the supports of two measured laminations $\lambda_1$ and 
$\lambda_2$ are not homotopic in $M$.
We shall explain how to deal with quasi-Fuchsian groups corresponding to
non-orientable surfaces later in \S5.

Our argument to prove Theorem \ref{converge} is basically as follows.
First we shall apply  the Theorems \ref{Thurston} and \ref{KS} to subgroups of $G$
corresponding to each component of the characteristic compression body of $C$ and the
complement of the characteristic compression body.
Then we shall regard the group $G$ as
constructed by gluing these convergent groups, and show the entire group converges by
analysing the group corresponding to the gluing surfaces. 

Now, we start the proof. 
Let $C$ be a convex core of $\H^3/G$, which is also a compact core since $G$ is geometrically
finite and has no parabolic elements.
Because of the existence of the ``nearest point retraction'' from $M$ to $C$ (refer to Morgan
\cite{Mo}), there is a homeomorphism from $C$ to $M$ isotopic to the inclusion,
which is the homotopical inverse of the nearest point retraction.
We call this homeomorphism the {\em natural homeomorphism} from $C$ to $M$.
Let $V$ be the characteristic compression
body of
$C$, and 
$V_1,\ldots , V_n$ its components.
Let $G^{V_k}$ be a subgroup of $G$ corresponding to the image of 
$\pi_1(V_k)$ in $\pi_1(M) \cong G$.
The quasi-conformal deformation $(G_i, \phi_i)$ induces that of 
$G^{V_k}$, which we denote by $(G_i^{V_k}, \phi_i|G^{V_k})$, where 
$G_i^{V_k}$ is nothing other than $\phi_i(G^{V_k})$.

The first step of the proof is to show that $\{(G_i^{V_k}, 
\phi_i|G^{V_k})\}$ converges in $AH(G^{V_k})$.

\begin{lem}
\label{compression converge}
The sequence $\{(G_i^{V_k}, \phi_i|G^{V_k})\}$ converges in 
$AH(G^{V_k})$ after passing to a subsequence.
\end{lem}
\begin{proof}
The compression body $V_k$ lifts homeomorphically to a compact core 
$\tilde{V}_k$ of $\mathbf{H}^3/G^{V_k}$.
Since $G^{V_k}$ is geometrically finite and has no parabolic 
elements, the hyperbolic 3-manifold 
$\H^3/G^{V_k}$ is homeomorphic to the interior of  a compact 
$3$-manifold $V'_k$, which is obtained as a Kleinian
manifold $(\H^3
\cup \Omega_{G^{V_k}})/G^{V_k}$.
The natural homeomorphism maps $\tilde{V}_k$ homeomorphically to $V_k'$.
Let $S$ be the exterior boundary of $V_k'$.
Then we can identify $S$ with a boundary component $S_j$ of $M$ via the corresponding
surfaces on the boundaries of $\tilde{V}_k$ and $V_k$. 
Recall that in Theorem \ref{converge}, a marked conformal structure 
$g^i_j$ was given on $S_j$.
We regard this $g^i_j$ as a marked conformal structure on $S$ by the identification given
above.
Then the quasi-conformal deformation $(G_i^{V_k},\phi_i|G^{V_k})$ is the image 
of $(g_j^i, h_i) \in \mathcal{T}(S) \times \mathcal{T}(\partial_i 
V'_k)$ by the Ahlfors-Bers map for some $h_i \in 
\T(\partial_i V'_k)$.

By assumption, $\{g_j^i\}$ either converges in the Thurston compactification of the Teichm\"{u}ller space to a maximal and 
connected projective 
lamination, which is contained in the projectivised Masur domain, or 
stays in a compact set of the Teichm\"{u}ller space.
In the former case, we can see that $\{(G_i^{V_k}, \phi_i|G^{V_k})\}$ 
converges after taking a subsequence by Theorem \ref{KS}.
In the latter case, by Theorem 2.1  in Canary \cite{CaD}, for every 
element $\gamma \in G^{V_k}$ represented by a closed 
curve on $S$, the 
translation length of $\phi_i(\gamma)$ is bounded as $i 
\rightarrow \infty$.
Since $G^{V_k} \cong \pi_1(V'_k)$ is carried by $\pi_1(S)$, this implies that 
$\{(G_i^{V_k},\phi_i|G^{V_k})\}$ converges after passing to a subsequence. 
\end{proof}

Let $(\Gamma^{V_k},\psi^{V_k})$ be the algebraic limit of a convergent subsequence of
\linebreak
$\{(G_i^{V_k}, \phi_i|G^{V_k})\}$.
Recall that the closure of the complement of $V$ in $C$ is a boundary-irreducible 3-manifold,
which we shall denote by $W$. 

We note that the compact core $C$ of $\mathbf{H}^3/G$ can be 
constructed from the components of $W$ and the components of $V$ by 
pasting them along incompressible boundary components.
If we allow ourselves to paste two boundary components of  $V$ each other, then
we can ignore the components of $W$ which are product $I$-bundles over 
closed surfaces.
For, the homeomorphism type does not change by gluing 
such a component one of whose boundary components does not lie on $\partial_i V$,
and
we can realise the same manifold by pasting two boundary components of $V$ if
such a component has both boundary components on $\partial_i V$.

Let $W_1, \ldots, W_m$ be the components of $W$ which are not
product 
$I$-bundles over closed surfaces.
For each  $W_j$, we define $G^{W_j}$ to be a subgroup of $G$ 
corresponding to the image of $\pi_1(W_j)$ in $\pi_1(C) \cong G$.
Then $W_j$ can be lifted homeomorphically to a compact core $\tilde{W}_j$ 
of $\mathbf{H}^3/G^{W_j}$.
Let $(G_i^{W_j},\phi_i|G^{W_j})$ be the quasi-conformal deformation of $G^{W_j}$ 
induced from $(G_i,\phi_i)$.
We denote by $\Phi_i^{W_j}$ a homotopy equivalence from $\H^3/G^{W_i}$ to
$\H^3/G_i^{W_j}$ corresponding to
$\phi_i|G^{W_j}$.

\begin{lem}
\label{irreducible converge}
The sequence $\{(G_i^{W_j}, \phi_i|G^{W_j})\}$ converges in 
$AH(G^{W_j})$ after passing to a subsequence.
\end{lem}

\begin{proof}
Suppose that $\{(G_i^{W_j}, \phi_i|G^{W_j})\}$ does not have a 
convergent subsequence.
Then, there are a boundary component 
$\tilde{T}$ of $\tilde{W}_j$ and disjoint unions of essential annuli 
$A^1_i, \dots , A^l_i$ 
properly embedded in $\tilde{W}_j$ with weights $w^1_i, \ldots , 
w^l_i$ for each $i$ as in Theorem \ref{Thurston},
such that $(w_i^1 A_i^1 \sqcup \dots \sqcup w_i^lA_i^l) \cap \tilde{T}$, regarded as
an element of $\ML(\tilde{T})$, converges to a (non-empty) measured 
lamination $\mu$
on $\tilde{T}$ as $i \rightarrow \infty$.
(In other words, we are considering $\tilde{T}$ which intersects the
codimension-1 lamination $L$ which appeared in the proof of Theorem
\ref{Thurston}.
The measured lamination $\mu$ is nothing but $L \cap \tilde{T}$.)
Let $T$ be the image of $\tilde{T}$ in $W_j$ under the covering projection.
The surface $T$ is either contained in $\partial C$ or an interior 
boundary component of some component $V_k$ of $V$.

First consider the case when $T$ is contained in $\partial C$.
Recall that in the situation of Theorem \ref{converge}, we are 
given a sequence of marked conformal structures $\{g_i\}$ 
on $T$ if we identify $C$ with $M$ by the natural homeomorphism.
 This $g_i$ is regarded as a point in $\T(\tilde{T})$ by identifying $\tilde{T}$ with $T$.
By assumption,  $\{g_i\}$ either stays in a compact set of $\T(\tilde{T})$ or converges  in
the Thurston compactification to a maximal and connected projective 
lamination $[\lambda]$ on $\tilde{T}$.
Suppose first that the former is the case.
Let $c$ be an essential simple closed curve on $\tilde{T}$ such that
$i(c,\mu)>0$.
The length of the closed geodesic freely homotopic to $\Phi_i^{W_j}(c)$ is 
bounded as $i \rightarrow \infty$ by Sullivan's theorem (see 
Epstein-Marden \cite{EM} for the proof) since the geodesic length of $c$ with 
respect to $g_i$ is bounded as $i \rightarrow \infty$.
This contradicts the last sentence of the conclusion of Theorem \ref{Thurston} 
asserting that the length of  such a closed geodesic  goes to 
infinity as $i \rightarrow \infty$.
 
Next suppose that the latter is the case.
Then, by Theorem \ref{length to 0}, there exist a sequence of 
simple closed curves $\{c_i\}$ and positive real numbers $\{r_i\}$ 
such that $\{r_i c_i\}$ converges to $\lambda$ and 
$r_i\mathrm{length}_{g_i}(c_i) \rightarrow 0$.
Let $c_i^*$ be the closed geodesic in $\H^3/G^{W_j}_i$ freely homotopic to
$\Phi_i^{W_j}(c_i)$.
Again by Sullivan's theorem proved by Epstein-Marden \cite{EM}, this implies that $r_i 
\mathrm{length}(c_i^*) \rightarrow 0$.
By applying Theorem \ref{Thurston} to $G^{W_j}$, we see that this is possible only when 
$i(\mu,\lambda)=0$.
 Since $\lambda$ is assumed to be maximal and 
connected,   the supports of $\lambda$ and 
$\mu$ coincide; hence in particular $\mu$ is also maximal and connected.
 (Here the assumption that $\lambda$ is maximal and connected is crucial.)
Since $\mu$ is the limit of $(w_i^1 A^1_i \sqcup \ldots \sqcup w_i^l A_i^l) \cap 
\tilde{T}$, this happens only when 
the characteristic pair of $\tilde{W}_j$ is an $I$-bundle over a 
closed surface, coinciding 
with the entire $\tilde{W}_j$.
Since we are considering the case when $T$ is a boundary component 
of $C$, the component $W_j \cong \tilde{W}_j$ can be a twisted $I$-bundle only when 
$C=W_j$. 
This is the case when $G$ is a quasi-Fuchsian group corresponding to a
non-orientable surface, which we are excluding by assumption now.
As we assumed that $W_j$ is not a product 
$I$-bundle over a closed surface, $\tilde{W}_j$ cannot be a product 
$I$-bundle either.
 Thus, in either case, we are lead to a contradiction.

Next we consider the case when $T$ is an interior boundary component 
of a component $V_k$ of $V$.
 Since $\tilde{T}$ intersects the limit $L$ of $w_i^1 A_i^1 \sqcup \dots
\sqcup
w_i^k A_i^k$ essentially, it contains a simple closed curve $c$ intersecting $L$
essentially.
The translation length of $\phi_i(c)$ goes to infinity by Theorem
\ref{Thurston}.
Then we see that for a subgroup
$G^T$ of $G^{W_j}$ associated to the image of $\pi_1(T)$ in 
$\pi_1(W_j) \subset \pi_1(C)$, the sequence $\{(\phi_i(G^T), \phi_i|G^T)\}$ cannot converge as
$i
\rightarrow 
\infty$ even if we take  a subsequence.
Since there is a conjugate of $\phi_i(G^T)$ contained in $G_i^{V_k}$, this 
contradicts Lemma \ref{compression converge}.
\end{proof}

Let $C_1$ and $C_2$ be two compact submanifolds of $C$, whose 
inclusions induce monomorphisms between the fundamental groups, such that 
$C_1 \cap C_2$ is a connected closed surface $T$ which is incompressible in $C$.
Let  $G^1$ be a subgroup of $G$ corresponding to the image of 
$\pi_1(C_1)$ in $\pi_1(C) \cong G$, and $G^2$ one corresponding to the image of $\pi_1(C_2)$.
Taking a conjugate of $G^2$, we can assume that $G^1 \cap G^2$ is a 
subgroup of $G$ corresponding to the image of $\pi_1(T)$ in 
$\pi_1(C) \cong G$.
Let $G'$ be a subgroup of $G$ obtained as the amalgamated free 
product of $G^1$ and $G^2$ over $G^1 \cap G^2$.
Suppose that the two sequences of quasi-conformal deformations $\{(\phi_i(G^1), 
\phi_i|G^1)\}$ and $\{(\phi_i(G^2), \phi_i|G^2)\}$ converge 
algebraically after passing to subsequences.

\begin{lem}
\label{composed convergence}
In this situation, the sequence of quasi-conformal deformations $\{ (\phi_i(G'), 
\phi_i|G') \}$ also converges in $AH(G')$ after passing to a subsequence.
\end{lem}

\begin{proof}
Since we have only to prove the existence of convergent subsequence, 
we shall take a subsequence each time it is necessary in the proof without mentioning it.
Since  $\{(\phi_i(G^1), \phi_i|G^1)\}$ converges in 
$AH(G^1)$, by taking conjugates, we can 
assume that the sequence of representations $\{\phi_i|G^1\}$ converges.
It is sufficient to prove that  $\{\phi_i|G'\}$ also converges then.

Since  $\{(\phi_i(G^2), \phi_i|G^2)\}$ also
converges, there are elements $t_i \in PSL_2\mathbf{C}$ such that 
$\{t_i(\phi_i|G^2) t_i^{-1}\}$ converges.
As $G^1 \cap G^2$ is isomorphic to $\pi_1(T)$ and $\{\phi_i|(G^1 \cap 
G^2)\}$ 
converges to a Kleinian group, we can choose 
elements $\gamma,\delta \in G^1 \cap G^2$  such that both 
$\{\phi_i(\gamma)\}$  and $\{\phi_i(\delta)\}$
converge to loxodromic elements, which do not commute each other.
We should note here that $\{t_i\phi_i(\gamma)t_i^{-1}\}$ and 
$\{t_i\phi_i(\delta)t_i^{-1}\}$ also 
converge to loxodromic elements of $PSL_2\mathbf{C}$.
This is possible only when $\{t_i\}$ converges to an element of 
$PSL_2\mathbf{C}$, as we can easily see the moves  under $t_i$ of the fixed points on $S^2_\infty$ 
of $\phi_i(\gamma)$ and $\phi_i(\delta)$.
Therefore, we see that $\{\phi_i|G^2\}$ also converges as 
representations.
Since $G'$ is generated by $G^1$ and $G^2$, it follows that 
$\{\phi_i|G'\}$ also converges.
\end{proof}

Next consider a submanifold $C_0$ in $C$ whose inclusion induces a 
monomorphism between the fundamental groups, and two of whose boundary 
components are parallel in $C$.
Let $C'$ be the submanifold of $C$ obtained by pasting the two 
parallel boundary components of $C_0$.
Let $G^{C_0}$ be a subgroup of $G$ corresponding to the image of 
$\pi_1(C_0)$ in $\pi_1(C) \cong G$.
Then a subgroup $G^{C'}$ corresponding to the image of $\pi_1(C')$ 
in $\pi_1(C) \cong G$ is obtained as an HNN-extension of $G^{C_0}$ 
over a subgroup corresponding to the fundamental group of one of the 
parallel boundary components of $C_0$.
Consider quasi-conformal deformations $\{(\phi_i(G^{C'}), 
\phi_i|G^{C'})\}$ of $G^{C'}$.

\begin{lem}
Suppose that the sequence $\{(\phi_i(G^{C_0}), \phi_i|G^{C_0})\}$ 
converges in \linebreak$AH(G^{C_0})$.
Then the sequence  $\{(\phi_i(G^{C'}) , 
\phi_i|G^{C'})\}$ also converges in $AH(G^{C'})$ after passing to a subsequence.
\end{lem}

\begin{proof}
Let $T_1$ and $T_2$ be boundary components of $C_0$ which are parallel in $C$ and pasted 
each other in $C'$.
We regard the fundamental groups $\pi_1(T_1)$ and $\pi_1(T_2)$ as contained 
in $\pi_1(C_0)$ by connecting a basepoint in $C_0$ to $T_1, T_2$ 
by arcs in $C_0$.
Then, the fundamental group $\pi_1(C')$ is generated by 
$\pi_1(C_0)$ and an element $t \in \pi_1(C')$ such that the 
conjugation by $t$ gives an isomorphism from $\pi_1(T_1)$ to 
$\pi_1(T_2)$ regarded as contained in $\pi_1(C')$.
Let $G^{T_1}$ and $G^{T_2}$ be the  subgroups of $G^{C_0}$ corresponding 
to $\pi_1(T_1)$ and $\pi_1(T_2)$ respectively, 
and we use the same symbol $t$ to denote the element of $G^{C'}$ 
corresponding to $t \in \pi_1(C')$.
By assumption, 
$\{\phi_i|G^{C_0}\}$ can be made to converge by taking conjugates; hence in particular $\{\phi_i|G^{T_1}\}$ 
and $\{\phi_i|G^{T_2}\}$ converge.
Since $G^{T_1}$ contains two elements both mapped to  loxodromic elements, 
which do not commute each other, by the same argument as in Lemma 
\ref{composed convergence}, we see that  $\{\phi_i(t)\}$ must also converge.
Thus, we have proved that $\{\phi_i|G^{C_0}\}$ also converges.
\end{proof}

Hence, if neither $V$ nor $W$ is empty, starting from quasi-conformal deformations of
subgroups of $G$ corresponding  to  components $V$ and non-product-$I$-bundle components of
$W$, for which  the convergence of the corresponding subgroups was proved in 
Lemmata \ref{compression converge}, \ref{irreducible converge}, and 
using the argument above repeatedly, we see that $\{(G_i, \phi_i)\}$ converges 
in $AH(G)$ after taking a subsequence.
For the case when $V$ or $W$ is empty, Theorem \ref{KS} above and 
Theorem 2.1 in \cite{OhI} imply Theorem \ref{converge}.
Thus, we have completed the proof of Theorem \ref{converge}.

\section{Proof of the main theorem}
For a geometrically finite group $G$ without parabolic elements and measured laminations $\lambda_{j_1}, \ldots
,
\lambda_{j_p}$ as given in Theorem \ref{main}, we construct a sequence 
of quasi-conformal deformations $\{(G_i, \phi_i)\}$ of $G$ as follows.
Let $S_1, \ldots , S_m$ be the boundary components of the 
compactification $M$ of $\H^3/G$, as in the statement of the theorem.
For each $S_{j_k}$ among the boundary components $S_{j_1}, \ldots , S_{j_p}$ on which measured 
laminations $\lambda_{j_k}$ are given, consider a sequence of marked conformal structures 
$\{g^i_{j_k}\}$ on $S_{j_k}$ which 
converges to the projective lamination $[\lambda_{j_k}]$ in the 
Thurston compactification of the Teichm\"{u}ller space.
For the remaining boundary components $S_{i_1}, \ldots S_{i_q}$, we 
define marked conformal structures 
$g^i_{i_1}, \ldots g^i_{i_q}$ to be the given marked conformal 
structures $m_1, \ldots , m_q$ in the statement, 
which are constant with respect to $i$.
Let $q : \mathcal{T}(\partial M) \rightarrow QH(G)$ be the Ahlfors-Bers
map.
We define a quasi-conformal deformation $(G_i, \phi_i)$ to be 
$q(g^i_1, \ldots , g^i_m)$ for the marked conformal structures defined above.

By Theorem \ref{converge}, this sequence $\{(G_i, \phi_i)\}$ 
converges to a Kleinian group $(\Gamma, \psi)$ in $AH(G)$ passing to 
a subsequence provided that $G_i$ is not quasi-Fuchsian.
When $G_i$ is quasi-Fuchsian, if the corresponding closed surface is
orientable, we
can obtain the same result making
 use of Theorem 2.4 in \cite{OhI}, instead of Theorem \ref{converge}.
 (Actually this can be proved by the argument of Thurston's
double limit theorem \cite{Th2}.)

 Suppose that $G$ is a quasi-Fuchsian group isomorphic to $\pi_1(S')$ for a
non-orientable surface $S'$.
Then $C$ is a twisted $I$-bundle over $S'$.
By Theorem \ref{Thurston}, there exists a sequence of
weighted unions of essential annuli $w_i^1 A_i^1 \sqcup \dots \sqcup w_i^kA_i^k$
for which the conclusion of the theorem holds; i.e., the limit of $w_i^1 A_i^1 \sqcup \dots \sqcup w_i^kA_i^k$ describes the behaviour of the three-dimensional geodesic lengths of closed curves on the boundary.
Since $\partial C$ is identified with $S_{j_1}$, consider the limit of  $S_{j_1}
\cap (w_i^1 A_i^1 \sqcup \dots \sqcup w_i^kA_i^k)$ in the measured
lamination space and denote it by $\mu$.
Then obviously $\mu$ is a lift of a measured lamination which represents
$S'
\cap (w_i^1 A_i^1 \sqcup \dots \sqcup w_i^k A_i^k)$ if we identify $S'$ with
its cross section embedded in $C$.
Since $\lambda_{j_1}$ is maximal and connected, it has non-trivial intersection
with any measured lamination except for those with the same supports as
$\lambda_{j_1}$.
By assumption, $\lambda_{j_1}$ is not a lift of a measured lamination on
$S'$; hence its support is different from that of $\mu$.
Therefore we have $i(\lambda_{j_1}, \mu) >0$.
This implies that if a sequence of weighted simple closed curves $\{w_i
\gamma_i\}$ converges to $\lambda_{j_1}$, then we have
$w_i\mathrm{length}_{\H^3/G_i}(\Phi_i(\gamma_i)) \rightarrow \infty$ by
Theorem \ref{Thurston}.
By Sullivan's theorem proved in \cite{EM}, this leads to
$w_i\mathrm{length}_{g^i_{j_1}}(\gamma_i) \rightarrow \infty$, which
contradicts Theorem \ref{length to 0}.
Thus we see that $\{(G_i, \phi_i)\}$ converges also in this case.
 
Thus, in every case, we have a limit $(\Gamma,\psi)$ of $\{(G_i,\phi_i)\}$
after taking a subsequence.
Let $C$ be the convex core of $\H^3/G$ as before.
Let $\Psi : \mathbf{H}^3/G \rightarrow \mathbf{H}^3/\Gamma$ be a 
homotopy equivalence inducing the isomorphism $\psi$ between 
$\pi_1(\mathbf{H}^3/G) \cong G$ and $\pi_1(\mathbf{H}^3/\Gamma) \cong 
\Gamma$, following our convention of the notation.
To complete the proof of Theorem \ref{main}, it is sufficient to prove the following
proposition.

\begin{prop}
\label{tame}
Regard $S_1, \ldots, S_m$ also as the boundary components of $C$ identifying $C$ with $M$ by
the natural homeomorphism, and let
$C'$ be a compact core of $\H^3/\Gamma$.
Then the restriction of the homotopy equivalence $\Psi|C$ can be homotoped 
to a homeomorphism $h: C \rightarrow C'$, and the following hold.
\begin{enumerate}
\item For each $k=1, \ldots, p$, the lamination $h(\lambda_{j_k})$ is unrealisable in $\H^3/\Gamma$.
Hence in particular $\Gamma$ is geometrically infinite.

\item Unless $G$ is free, $\Gamma$ has no parabolic elements and $T_{j_k}=h(S_{j_k})$ faces a geometrically infinite and topologically tame end.
\item
For each of the remaining components, $S_{i_k}$, $k=1,\ldots, q$, the 
end facing $h(S_{i_k})$ 
is geometrically finite.

\item $\Gamma$ is topologically tame.
\end{enumerate}

Therefore, in particular, $\H^3/\Gamma$ is compactified to a 
3-manifold $M'$, which is homeomorphic to $C'$.
The homeomorphism $h$ extends to a homeomorphism $\overline{h}: M 
\rightarrow M'$, and the corresponding conformal structure at 
infinity on $\overline{h}(S_{i_k})$ is $\overline{h}_*(m_k)$.
\end{prop}

\begin{proof}[Proof of Proposition \ref{tame}]
Let $c_i$ be the shortest simple closed geodesic with respect to the 
hyperbolic structure on $(S_{j_k}, g^i_{j_k})$.
By the compactness of \linebreak$\PL(S)$, passing to a subsequence, there is a sequence of positive real 
numbers $r_i$ such that $\{r_i c_i\}$ converges to a measured 
lamination $\mu$.
Obviously we can assume that either $r_i \rightarrow 0$, or $r_i$ is 
constant and $\mu$ is a simple closed curve, by taking a subsequence.
Since $c_i$ is the shortest simple closed geodesic, 
$\l_{g^i_{j_k}}(c_i)$ is bounded above independently of $i$.
If $i(\lambda_{j_k}, \mu) >0$, then by a well-known fact about the 
Thurston compactification (see expos\'{e} 8 of \cite{FLP} for instance), we have $r_i
\l_{g_{j_k}^i}(c_i) \rightarrow 
\infty$, and we get a contradiction.
Therefore $i(\lambda_{j_k}, \mu)=0$, which implies the supports of 
$\lambda_{j_k}$ and $\mu$ coincide by the assumption that 
$\lambda_{j_k}$ is maximal and connected.
In particular we see that $\mu$ is not a simple closed curve, hence $r_i 
\rightarrow 0$, and that $\mu$ is contained in $\M(S_{j_k})$.

On the other hand, since the length of $c_i$ with respect to 
$g^i_{j_k}$ is bounded above as $i \rightarrow \infty$, its extremal length 
is also bounded above, and by Corollary 2 in Maskit \cite{MaF}, we see that 
$\{\l_{\H^3/G_i}(\Phi_i(c_i))\}$ is bounded; hence 
$r_i\l_{\H^3/G_i}(\Phi_i(c_i)) \rightarrow 0$ as $i \rightarrow 
\infty$.
By Theorem \ref{continuity},
it follows 
that $\mu$ cannot be realised by a pleated surface homotopic 
to $\Psi|S_{j_k}$.
(Recall that $\mu$ is maximal and connected.)
Since $\lambda_{j_k}$ has the same support as $\mu$, the measured 
lamination $\lambda_{j_k}$ cannot be realised by such a pleated surface 
either.
Therefore once we can show $\Phi|C$ is homotoped to a homeomorphism,  (1) will follow.

Now assume that $G$ is not free.
Consider a subgroup $\Gamma^{S_{j_k}}$ of $\Gamma$ corresponding to 
$(\Psi|S_{j_k})_\# \pi_1(S_{j_k})$ in  $\pi_1(\mathbf{H}^3/\Gamma) \cong 
\Gamma$, where $\#$ denotes the induced homomorphism.
Let $C'_{S_{j_k}}$ be a compact core of $\mathbf{H}^3/\Gamma^{S_{j_k}}$.
The map $\Psi|S_{j_k}$ can be lifted to a map $\tilde{f} : S_{j_k} 
\rightarrow \H^3/\Gamma^{S_{j_k}}$.
Since the measured lamination $\lambda_{j_k}$ cannot be realised by a 
pleated surface homotopic to $\Psi|S_{j_k}$, it cannot be realised by 
a pleated surface homotopic to $\tilde{f}$ either.
Therefore, as was shown in \cite{Ot}, there are pleated 
surfaces $\tilde{f}_n : S_{j_k} \rightarrow \mathbf{H}^3/\Gamma^{S_{j_k}}$ 
homotopic to $\tilde{f}$, which realise weighted simple closed curves in the 
Masur domain converging to $\lambda_{j_k}$, and tend to an end $\tilde{e}$ as $n \rightarrow \infty$.
By Theorem 4.1  in \cite{OhM}, this implies that the end $\tilde{e}$ is geometrically infinite and has a neighbourhood $\tilde{E}$ homeomorphic to $\Sigma \times 
\mathbf{R}$ such that $\tilde{f}$ is homotopic (in $\H^3/\Gamma^{S_{j_k}}$) to a
covering map 
onto $\Sigma \times \{t\}$.
By cutting $\mathbf{H}^3/\Gamma^{S_{j_k}}$ along $\Sigma \times 
\{t\}$, and applying the relative core theorem due to McCullough 
\cite{Mc}, we can see that there 
is a compact core $K$ containing $\Sigma \times \{t\}$ as a boundary 
component.
The uniqueness of core  implies that there is a 
homeomorphism between the cores $K$ 
and $C'_{S_{j_k}}$, which induces an inner-automorphism of 
$\pi_1(\mathbf{H}^3/\Gamma^{S_{j_k}})$.
In particular, there must be a boundary component $T$ of $C'_{S_{j_k}}$ 
which is homeomorphic to $\Sigma$, and is homotopic  to $\Sigma\times 
\{t\}$  in $\mathbf{H}^3/\Gamma^{S_{j_k}}$.
Then, since the fundamental group of $T$ carries a conjugate of the image 
of $\pi_1(\tilde{f}(S_{j_k}))$, which must be equal to the entire 
$\pi_1(C'_{S_{j_k}})$, this can occur only when $C'_{S_{j_k}}$ 
is a compression body, and $T$ is homeomorphic to $S_{j_k}$.
Thus, we see that $\tilde{E}$ is homeomorphic to $S_{j_k} \times 
\mathbf{R}$ in such a way that $\tilde{f}(S_{j_k})$ is homotopic in
$\H^3/\Gamma^{S_{j_k}}$ to
the 
surface corresponding to $S_{j_k} 
\times \{t\}$.

Let $p : \mathbf{H}^3/\Gamma^{S_{j_k}} \rightarrow 
\mathbf{H}^3/\Gamma$ be the covering map associated to the inclusion 
$\Gamma^{S_{j_k}} \subset \Gamma$.
The covering theorem of Canary \cite{CaT}, which is a generalisation 
of Thurston's covering theorem, implies that we can take a 
neighbourhood $\tilde{E}$ of $\tilde{e}$ homeomorphic to $S_{j_k} 
\times \mathbf{R}$ as above so that $p|\tilde{E}$ is a finite-sheeted 
covering to its image; a neighbourhood $E$ of an end $e$ of 
$\mathbf{H}^3/\Gamma$, which is homeomorphic to $\overline{\Sigma} 
\times \mathbf{R}$ for a closed surface $\overline{\Sigma}$ in such a way that 
$p|S \times \{t\}$ covers $\overline{\Sigma} \times \{t\}$.
By the same argument as above, we can see that there is a boundary 
component of the compact core $C'$, which is homotopic to 
$\overline{\Sigma} \times \{t\}$.
This also implies that the end $e$ is geometrically infinite, topologically tame, and contains no cusps.
This will show (2) once $\Phi|C$ is shown to be homotopic to a homeomorphism to $C'$.

Thus, unless $G$ is free, we have shown that for each $S_{j_k}$ among $S_{j_1}, \ldots, 
S_{j_p}$, the homotopy equivalence $\Psi$ can be homotoped so that 
$\Psi|S_{j_k}$ is a finite-sheeted covering onto a boundary 
component of the compact core $C'$.

Next we consider a boundary component $S_{i_k}$ on which we fixed a 
conformal structure to construct $(G_i, \phi_i)$.
As before, we let $G^{S_{i_k}}$ be a subgroup of $G$ corresponding to 
the image of $\pi_1(S_{i_k})$ in $\pi_1(\mathbf{H}^3/G) \cong G$.
There is a unique component $\overline{S}_{i_k}$ of $\Omega_G/G$ which 
lies at infinity of the end facing $S_{i_k}$ in the Kleinian manifold 
$(\H^3\cup \Omega_G)/G$.
Let $\Omega_0$ be a component of $\Omega_G$ invariant by 
$G^{S_{i_k}}$, which is a lift of $\overline{S}_{i_k}$ in $\Omega_G$.
Let $f_i : S^2_\infty \rightarrow S^2_\infty$ be a quasi-conformal 
homeomorphism inducing the quasi-conformal deformation $(G_i,\phi_i)$.
Then $f_i|\Omega_0$ is conformal for all $i$.
As was shown in Lemma 3 in Abikoff \cite{Ab}, this implies that there 
is a component $\Omega'_0$ of the region of discontinuity 
$\Omega_\Gamma$, which is invariant by $\psi(G^{S_{i_k}})$, such that 
$\Omega_0/G$ is conformal to $\Omega_0'/\Gamma$ preserving the 
markings.
(Here the markings make sense only up to the ambiguity caused by
auto-homeomorphisms  of the boundaries of Kleinian manifolds which are
homotopic to the identity
in the Kleinian manifolds.)
Let $S'_{i_k}$ be a boundary component of $C'$ facing the end which
corresponds to $\Omega'_0/\psi(G^{S_{i_k}})$.
Then $S_{i_k}'$ faces a geometrically finite end, whose corresponding 
marked conformal structure at infinity is $m_k$ with the marking 
determined by $\Psi|S_{i_k}$, and $\Psi|S_{i_k}$ is homotopic to a 
homeomorphism onto $S'_{i_k}$ via the natural homeomorphism and the nearest point retraction.

Combining this with the preceding argument, we see that we can 
homotope $\Psi$ so that each boundary component of $C$ is mapped to a 
boundary component of $C'$ unless $G$ is free.
As $C'$ is a compact core, evidently we can assume that $\Psi(C) 
\subset C'$.
Since $\Psi|C$ is a homotopy equivalence to $C'$, by Waldhausen's 
theorem \cite{Wa}, we can further homotope $\Psi$ to a homeomorphism 
$h$ from $C$ to $C'$ unless $G$ is free.
In the case when $G$ is free, both $C$ and $C'$ must be handlebodies and any homotopy equivalence from $C$ to $C'$ can be homotoped to a homeomorphism.
Thus we have completed the proof of (1)-(3).

It only remains to show the topological tameness of $\Gamma$.
This has already been proved for the case when $G$ is not free by showing that every geometrically infinite end corresponding to $S_{j_k}$ is topologically tame and the remaining ends are geometrically finite.
In the case when $G$ is free, we need to invoke either the work of Brock-Souto \cite{BS} and Brock-Bromberg-Evans-Souto \cite{BBES} or a general resolution of Marden's conjecture by Agol \cite{Ag} and Calegari-Gabai \cite{CG}.
\end{proof}

\section{Appendix}
We shall give a sketch of the proof of Theorem \ref{continuity}.
Since this theorem is derived directly from the continuity of length function
by Brock
\cite{Br} in the case when $S$ is incompressible, we have only to deal with the
case when $S$ is
compressible. The main ingredients for this argument are contained in Ohshika \cite{OhM}.

Suppose, seeking a contradiction, that there is a pleated surface homotopic to $\Psi|S$
realising $\mu$.
(Recall that we assumed that $\mu$ is maximal and connected.)
Then, as was shown in Lemma 4.7 in \cite{OhM}, for any $\epsilon >0$, there are a train track
$\tau$ with a weight system $\omega$ carrying $\mu$  and a map $f : S \rightarrow \H^3/\Gamma$
homotopic to
$\Psi|S$ which is adapted to $\tau$ such that the total curvature and the total quadratic
variation of angles of
$f(\tau)$ are less than $\epsilon$. Here a map $f$ is said to be adapted to a train track
$\tau$ when each branch of
$\tau$ is mapped to a geodesic arc and there is a tied neighbourhood $N_\tau$ of $\tau$ such
that each tie is mapped to a point by
$f$.
The total curvature of $f(\tau)$ is defined to be the sum of the exterior angles formed by the
images of two adjacent branches multiplied by the weight flowing from the first branch to the
second.
Similarly, the total quadratic variation of angles is the weighted sum of the squares of such 
exterior angles.

By taking a subsequence, we can assume that $\{G_i\}$ converges geometrically to a Kleinian
group $G_\infty$ containing $\Gamma$ as a subgroup.
Then, there is a  $(K_i, r_i)$-approximate isometry $\rho_i: B_{r_i}(\H^3/G_i, x_i)
\rightarrow B_{K_ir_i}(\H^3/G_\infty, x_\infty)$, where $x_i$ and $x_\infty$ are
the images of some basepoint $\tilde{x} \in \H^3$, such that
$K_i^{-1}d(x,y)
\leq \linebreak d(\rho_i(x),
\rho_i(y))
\leq K_id(x,y)$ for every $x, y \in   B_{r_i}(\H^3/G_i, x_i)$, with $r_i \rightarrow \infty$
and $K_i \rightarrow 1$.
It is easy to see that for sufficiently large $i$, we can construct a map $f_i$ homotopic to
$\rho_i^{-1}
\circ f$ which is adapted to $\tau$ by defining for each branch $b$ of $\tau$, its image
$f_i(b)$ to be the geodesic arc homotopic to
$\rho_i^{-1}f(b)$ fixing the endpoints.
Also, we can easily show that the total curvature and the total quadratic variation of angles
for $f_i(\tau)$ are less than $2\epsilon$ for sufficiently large $i$.
(See the proof of Lemma 6.10 in \cite{OhM}.)

Now, as we assumed that $\mu$ is maximal and connected, we can take a train
track $\tau$ as above so that every measured lamination near $\mu$ is carried
by $\tau$.
Then by the same argument as a proof of   Proposition 5.1 of Bonahon
\cite{BoA} (Lemma 6.10 in
\cite{OhM}),
we see that for every $\epsilon>0$, there are $i_0$ and a neighbourhood $U$
of
$\mu$ in
$\mathcal{M}(S)$ such that for any weighted simple closed curve $\gamma \in U$, the closed
geodesic $\gamma_i^*$ homotopic to $\Phi_i(\gamma)$ has a part with length
$(1-\epsilon)\l(\gamma_i^*)$ lying in the $\epsilon$-neighbourhood of
$f_i(\tau)$.
By applying the same lemma fixing $i$ and considering a realisation of
$\lambda_i$ by a pleated surface homotopic to $\Phi_i|S$, we see that this gives a positive
lower bound for the length of the realisation of $\lambda_i$ in $\H^3/G_i$.
This implies that
$\ell(\lambda_i, G_i)$ is bounded below by a positive constant as $i \rightarrow \infty$,
contradicting our assumption.


\begin{thebibliography}{99}
\bibitem{Ab}
W.\ Abikoff, On boundaries of Teichm\"{u}ller spaces and on Kleinian 
groups, Acta Math.\ {\bf 134}, (1975), 211--237.  
\bibitem{Ag}
I.\ Agol, Tameness of hyperbolic 3-manifolds, arXiv:math.GT/0405568
 \bibitem{AC}
 J.\ W.\ Anderson and R.\ D.\ Canary, Algebraic limits of Kleinian groups which rearrange the
pages of a book, Invent. Math. {\bf 126} (1996),  205--214.
\bibitem{Be}
M.\ Bestvina, Degenerations of the hyperbolic space. Duke Math. J. {\bf 56}
(1988), 143--161.
\bibitem{BoE}
F.\ Bonahon, Cobordism of automorphisms of surfaces, 
Ann.\ Sci.\ \'{E}cole Norm.\ Sup (4) {\bf 16}, (1983), 237-270.
\bibitem{BoA}
\bysame, Bouts des vari\'{e}t\'{e}s hyperboliques de dimension 3, Ann.\ of 
Math.\ {\bf 124}, (1986) 71-158.
\bibitem{Br}
J.\ Brock, Continuity of Thurston's length function, Geom.\ Funct.\ 
Anal.\ {\bf 10}, (2000),  741--797. 
\bibitem{BB}
J.\ Brock and K.\ Bromberg, On the density of geometrically finite Kleinian groups.  Acta Math.  {\bf 192}  (2004),  33--93.
\bibitem{BBES} J.\ Brock, K.\ Bromberg, R.\ Evans, and J.\ Souto, Tameness on the boundary and Ahlfors' measure conjecture,  Publ. I.H.E.S. {\bf 98} (2003), 145-166
\bibitem{BCM} J.\ Brock, R.\ Canary, and Y.\ Minsky, The classification of Kleinian surface groups, II: The Ending Lamination
Conjecture, preprint, arXiv math.GT/0412006.
\bibitem{BS}
J.\ Brock and J.\ Souto, Algebraic limits of geometrically finite manifolds are tame,
 Geom. Funct. Anal.  {\bf 16}  (2006), 1--39.
\bibitem{Brom} K.\ Bromberg, Hyperbolic cone-manifolds, short geodesics, and Schwarzian derivatives.  J. Amer. Math. Soc.  {\bf 17}  (2004), 783--826.
\bibitem{CG}D.\ Calegari and D.\ Gabai, Shrinkwrapping and the taming of hyperbolic 3-manifolds,  J. Amer. Math. Soc.  {\bf 19}  (2006),  385--446
\bibitem{CaD}
R.\ Canary, The Poincar\'{e} metric and a conformal version of a theorem 
of Thurston, Duke Math.\ J.\ {\bf 64} (1991), 349--359.
\bibitem{CaJ}
\bysame, Ends of hyperbolic $3$-manifolds, J.\ Amer.\ Math.\ Soc.\ 
{\bf 6}, (1993), 1--35.
\bibitem{CaT}
\bysame, A covering theorem for hyperbolic 3-manifold and its 
applications, Topology {\bf 35}, (1996), 751-778.
\bibitem{CM} R.\ Canary and D.\ McCullough, Homotopy equivalences of 3-manifolds and deformation theory of Kleinian groups,  Mem. Amer. Math. Soc.  {\bf 172}  (2004),  no. 812, xii+218 pp.
  \bibitem{CS} M.\ Culler\ and\ P.\ B.\ Shalen,  Varieties of group representations and
splittings of $3$-manifolds, Ann.\ of Math.\ (2) {\bf 117} (1983),  109--146.
\bibitem{EM} D.\ B.\ A.\ Epstein and A.\ Marden, Convex hulls in 
hyperbolic space, a theorem of Sullivan , and measured pleated 
surfaces, {\em Analytical and geometric aspects of hyperbolic spaces}, 
London Math.\ Soc.\ Lecture Note Ser.\ {\bf 111} Cambridge Univ.\ 
Press (1987), 113-253.
\bibitem{FLP} A.\ Fathi, F.\ Laudenbach, et V.\ Po\'{e}naru, Travaux 
de Thurston sur les surfaces, Ast\'{e}risque, {\bf 66-67}, (1979).
\bibitem{JS} W.\ Jaco and P.\ Shalen, Seifert fibered spaces in 
$3$-manifolds, Mem.\ Amer.\ Math.\ Soc.\ {\bf 21} (1979), no.\ 220.
\bibitem{Jo} K.\ Johannson, Homotopy equivalences of $3$-manifolds 
with boundaries, Lecture Notes in Mathematics, {\bf 761} Springer, 
Berlin, 1979.
\bibitem{KS} G.\ Kleineidam and J.\ Souto, Algebraic convergence of 
function groups, Comment. Math. Helv. {\bf 77} (2002),  244--269.
\bibitem{KS2}\bysame,\bysame,  Ending laminations in the Masur domain, {\it
Kleinian groups and hyperbolic
3-manifolds (Warwick, 2001)}, 105--129, Cambridge Univ. Press, Cambridge, 2003
\bibitem{Ma}
A.\ Marden, The geometry of finitely generated Kleinian groups, Ann.\ 
of Math.\ {\bf 99}, (1974), 465-496.
\bibitem{MaF}
B.\ Maskit, Comparison of hyperbolic and extremal lengths, Ann.\ Acad.\ 
Sci.\ Fenn.\ Ser.\ A I Math.\ {\bf 10} (1985), 381--386.
\bibitem{Mc}
D.\ McCullough, Compact submanifolds of 3-manifolds with 
boundary, Quart.\ J.\ Math.\ Oxford {\bf 37}, (1986), 299-307.
\bibitem{MMS}
D.\ McCullough, A.\ Miller and G.A.\ Swarup, Uniqueness 
of cores of non-compact 3-manifolds, J.\ London Math.\ Soc.\ {\bf 32}, (1985),
548-556.
\bibitem{Mi} Y.\ Minsky, The classification of Kleinian surface groups, I: Models and bounds, preprint, arXiv math.GT/0302208
  \bibitem{Mo} J.\ Morgan, On Thurston's uniformisation theorem for 
three-dimensional 
manifolds, {\em The Smith conjecture}  Academic Press (1984), 37--125.
\bibitem{MS} J.\ Morgan and P.\ Shalen, Valuations, trees, and degenerations of hyperbolic structures. I.  Ann. of Math. (2) 
{\bf 120 } (1984),   401--476.
\bibitem{MS3} \bysame, \bysame,  Degenerations of hyperbolic structures. III.
Actions of
$3$-manifold groups on trees and Thurston's compactness theorem. Ann. of Math. (2)
{\bf 127} (1988),  457--519
\bibitem{NS} H.\ Namazi and J.\ Souto,  Nonrealizability in handlebodies and ending                                 
laminations, \emph{in preparation.}
\bibitem{OhZ}
K.\ Ohshika, On limits of quasi-conformal deformations of Kleinian groups.  Math. Z.  {\bf 201}  (1989),   167--176.
\bibitem{OhI}
\bysame, Limits of geometrically tame Kleinian groups, Invent.\ Math.\ 
{\bf 99}, (1990),185--203. 
\bibitem{OhL}
\bysame, Ending laminations and boundaries for deformation spaces 
of Kleinian groups, J.\ London Math.\ Soc.\ {\bf 42}(1990), 111--121.
\bibitem{OhE}
\bysame, Geometrically finite Kleinian groups and parabolic elements, 
Proc.\ Edinburgh Math.\ Soc.\ {\bf 41} (1998), 141--159.
\bibitem{OhG}
\bysame, Divergent sequences of Kleinian groups, {\em The Epstein birthday 
schrift} Geom.\ Topol.\ Monogr, {\bf 1}, Geom.\
Topol.\ , Univ.\ Warwick, Coventry, (1998), 419--450.
\bibitem{OhM}
\bysame, Kleinian groups which are limits of geometrically finite 
groups, Mem.\ Amer.\ Math.\ Soc., {\bf 177}, (2005), no.\ 834
\bibitem{Ohpr}
\bysame, Realising end invariants by limits of minimally parabolic
geometrically finite groups, arXiv:math.GT/0504546
\bibitem{Ot}
J.-P.\ Otal, Courants g\'{e}od\'{e}siques et produits libres, 
Th\`{e}se, Universit\'{e} de Paris-Sud, Orsay
\bibitem{OtA} \bysame, Le th\'{e}or\`{e}me d'hyperbolisation pour les vari\'{e}t\'{e}s fibr\'{e}es de dimension 3.   Ast\'{e}risque  {\bf 235}  (1996), x+159 pp. 
\bibitem{Pa}
F.\ Paulin, Topologie de Gromov \'{e}quivariante, structures hyperboliques et arbres r\'{e}els. Invent. Math.
{\bf 94} (1988),  53--80
\bibitem{Sc}
G.\ P.\ Scott, Compact submanifolds of $3$-manifolds, J.\ London 
Math.\ Soc.
(2) {\bf 7}, (1973), 246--250.
\bibitem{Sk} R.\ Skora, Splittings of surfaces. J. Amer. Math. Soc. {\bf 9}
(1996),  605--616.
\bibitem{Th2} W.\ Thurston, Hyperbolic structures on 3-manifolds II : 
Surface groups and 3-manifolds which fiber over the circle, preprint, 
ArXiv math.GT/9801045
\bibitem{Th3} \bysame, Hyperbolic structures on 3-manifolds III:
Deformation of 3-manifolds with incompressible boundary, preprint, 
ArXiv math.GT/9801058
\bibitem{Thm} \bysame, Minimal stretch maps between hyperbolic 
surfaces, preprint, ArXiv Math.GT/9801039
\bibitem{Wa}
F.\ Waldhausen, On irreducible $3$-manifolds which are sufficiently 
large, Ann.\ of Math.\ (2) {\bf 87} (1968) 56--88.
\end{thebibliography}
\end{document}